\documentclass[12pt,a4paper,reqno]{amsart}
\usepackage{cite}
\usepackage{amsrefs, tabularx}
\usepackage{caption}
\usepackage{subcaption}
\usepackage{amsmath,amssymb,amsfonts,amsthm,epsfig,graphicx,color}
\allowdisplaybreaks
\usepackage{esint}
\usepackage[margin=0.75in]{geometry}
\usepackage{url}
\usepackage{hyperref}

\usepackage{mathrsfs}

\theoremstyle{definition}

\numberwithin{equation}{section}
\numberwithin{figure}{section}

\usepackage{mathtools}
\mathtoolsset{showonlyrefs}

\renewcommand{\le}{\leqslant}

\renewcommand{\ge}{\geqslant}

\def\R{\mathbb R}
\def\N{\mathbb N}

\renewcommand{\epsilon}{\varepsilon}

\allowdisplaybreaks
\begin{document}

\title[Description of an ecological niche]{Description of an ecological niche\\
for a mixed local/nonlocal dispersal:\\
an evolution equation and a new Neumann condition
arising from the
superposition\\
of Brownian and L\'evy processes}\thanks{
{\em Serena Dipierro}:
Department of Mathematics
and Statistics,
University of Western Australia,
35 Stirling Hwy, Crawley WA 6009, Australia.
{\tt serena.dipierro@uwa.edu.au}\\
{\em Enrico Valdinoci}:
Department of Mathematics
and Statistics,
University of Western Australia,
35 Stirling Hwy, Crawley WA 6009, Australia. {\tt enrico.valdinoci@uwa.edu.au}\\
It is a pleasure to thank Michael Small for his comments and
Gianmaria Verzini for very
interesting
conversations about the Skorokhod problem.\\
The authors are members of INdAM and AustMS and
are supported by the Australian Research Council
Discovery Project DP170104880 NEW ``Nonlocal Equations at Work''.
The first author is supported by
the Australian Research Council DECRA DE180100957
``PDEs, free boundaries and applications''. 
The second author is supported by
the Australian Laureate Fellowship
FL190100081
``Minimal surfaces, free boundaries and partial differential equations''.}

\begin{abstract}
We propose here a motivation
for a mixed local/nonlocal problem with a new type
of Neumann condition.

Our description is based on formal expansions and approximations.
In a nutshell, a biological species is supposed to diffuse
either by a random walk or by a jump process, according to
prescribed probabilities. If the process makes an individual
exit the niche, it must come to the niche right away,
by selecting the return point according to the underlying stochastic process.
More precisely,
if the random particle exits the domain,
it is forced to immediately reenter the domain,
and the new point in the domain is chosen randomly
by following a bouncing process with 
the same distribution as the original one.

By a suitable definition outside the niche, the density of the population
ends up solving a mixed local/nonlocal equation, in which
the dispersion is given by the superposition of the classical and the fractional
Laplacian. This density function satisfies two types of Neumann
conditions, namely the classical Neumann condition
on the boundary of the niche, and a nonlocal
Neumann condition
in the exterior of the niche. 
\end{abstract}

\author{Serena Dipierro}
\author{Enrico Valdinoci}


\maketitle

\section{Introduction}

The goal of this note is to provide an intuitive mathematical
explanation related to a recent model proposed in~\cite{VERO}
to describe the diffusion of a {\em biological population}
living in an {\em
ecological niche} and subject to {\em
both local and nonlocal dispersals}.
The model is motivated by the biologically relevant situation of a population
following long-jump foraging patterns alternated with
focused searching strategies at small scales. The exterior of the niche is accessible
by the population, but presents a hostile environment that forces the population to
an immediate return to the niche after a possible egression.
\medskip

Namely, we present a model of an evolution equation driven by a diffusive operator of mixed
local and nonlocal type (that is, the sum of a classical and a fractional Laplacians)
with a perspective proper for application in ecology. The solution
of this evolution equation represents the density of a biological population
living in an ecological niche, which in turn corresponds to the domain in which the equation takes place.
The mixed operator is the outcome of a superposition of a
classical (i.e. Brownian) and a long-range (i.e. L\'evy) stochastic processes.
In particular, we are interested in describing the situation in which
the population exits the biological niche and immediately comes back to it by following
the above mentioned stochastic process: this phenomenon naturally leads to the study of
a new type of mixed Neumann condition made of two separate prescriptions
(namely, a classical Neumann condition on the boundary of the given domain
and a nonlocal Neumann condition set on the exterior of the domain).
\medskip

More specifically,
the main mathematical framework
presented in~\cite{VERO} can be described by the
diffusive equation (endowed with external data)
\begin{equation}
\begin{split}&
\label{EQUAZIONE}\frac{
\partial U}{\partial t}=\alpha\Delta U -\beta(-\Delta)^sU
\qquad{\mbox{ in }}\; \Omega\times(0,+\infty),\\&
{\mbox{where }}\qquad
(-\Delta)^s U(x):=\int_{\R^n}\frac{U(x)
-U(y)}{|x-y|^{n+2s}}\,dy,\end{split}
\end{equation}
\begin{align}&
\label{NEUM-LOC}
\lim_{\vartheta\to0^+}\frac{U(x+\vartheta\nu,t)-U(x,t)}{\vartheta}
=:
\frac{\partial U}{\partial \nu}(x)=0\qquad
{\mbox{ for every }}\;x\in\partial\Omega,\\&
\label{NEUM-NLOC}
\int_\Omega \frac{U(x)-U(y)}{|x-y|^{n+2s}}\,dy
=0\qquad    {\mbox{ for every }}\;x\in\R^n \setminus \overline\Omega.
\end{align}
In this setting, $U=U(x,t)$ is a smooth function
depending on the space variables~$x\in\R^n$ and~$t\in[0,+\infty)$,
possibly endowed with some initial condition at time~$t=0$.
The diffusive parameters~$\alpha$ and~$\beta$
are positive constants, $\Omega$ is
a, say, bounded and smooth domain in~$\R^n$,
with external unit normal given by~$\nu$,
and~$s\in(0,1)$.
Also, the integral in~\eqref{EQUAZIONE}
is evaluated in the principal value sense
whenever needed to cope with the singularity produced by
the vanishing of the denominator.
\medskip

The biological interpretation of the setting
in~\eqref{EQUAZIONE}--\eqref{NEUM-NLOC}
arises from a biological population with density~$u=u(x,t)$
living in the niche described by the domain~$\Omega$:
then, the function~$U$ in~\eqref{EQUAZIONE}--\eqref{NEUM-NLOC}
is obtained from~$u$ in such a way that~$U$ coincides with~$u$
in~$\Omega$ and it is a convenient extension of~$u$ outside~$\Omega$.
Hence, the knowledge of the function~$U$
in~\eqref{EQUAZIONE}--\eqref{NEUM-NLOC} completely determines
the density of the biological population in the niche
(in addition, as pointed out in~\cite{MR4102340}
the setting in~\eqref{EQUAZIONE}--\eqref{NEUM-NLOC}
for~$U$
can be reformulated as a regional problem
in~$\Omega$ for an integrodifferential operator for~$u$
whose kernel has a logarithmic singularity along~$\partial\Omega$,
see also~\cite{AUD}).\medskip

The population diffuses according to two types
of dispersals, namely a classical one, related to the classical Laplacian,
and a nonlocal one, modeled on L\'evy flights
and encoded by the fractional Laplacian
in~\eqref{EQUAZIONE}. While the occurrence of nonlocal
diffusion in biological species is an extensively studied phenomenon
(see e.g.~\cite{NATU, NATU2})
and it is commonly known that captivity is lethal for some species
used to cruise broad open distances
at a time (see e.g.~\cite{POL, POL2, SHARK} -- suggesting
that for such species the classical Neumann condition
is somehow as lethal as a homogeneous Dirichlet condition),
a precise description on how this is compatible with the existence
of a prescribed niche for the population is a topic of
current investigation, and the contribution in~\cite{VERO}
was precisely to propose a clear mathematical framework
to settle a model in which both classical and anomalous diffusions take
place,
the case of purely nonlocal diffusion being firstly
introduced in~\cite{dirova}.
To make the model biologically
richer, in~\cite{VERO} the dispersal operator was also complemented
by a logistic equation taking into account the competition for
resources and also by an auxiliary possible pollination (or mating call)
term allowing a rise of the birth rate in view of some
interactions between individuals located at some range from each other.\medskip

The importance of methods from probability and statistics to address the question of animal foraging also emerged in the analysis of the best strategy for the search of randomly located targets, see~\cite{89897765NATU, VISWANATHAN2001, MR2609393, MR3403266}. With respect to this, we recall that the foraging efficiency problem presents multifaceted features, depending, for instance, on whether the target is stationary or in motion (and with which velocity with respect to the forager), on whether or not the forager is capable of spotting the target at a finite distance or sense targets while moving, on whether or not the forager eliminates the targets when it finds them. The answer to the optimal foraging strategy is very sensible to all these details, as well as to the density of the targets and, of course, to the notion of optimality and foraging efficiency chosen to quantify and compare the different possible strategies: a possible paradigm emerging to different analyses seems to be that L\'evy processes provide optimal search strategies in case of sparse, randomly distributed, immobile, revisitable targets in unbounded domains, since a long-jump searching strategy avoids oversampling (or at least allows to sample regions more quickly)
in these situations, and in several destructive situations (when the forager eliminates the target) the ``ballistic'' limit (corresponding to the case~$s\to0^+$ in our notation) happens to offer the highest foraging performances -- conversely, Brownian strategies outperform L\'evy ones in case of densely distributed targets, see e.g.~\cite{Klages R} for a thorough review of different possible scenarios.\medskip

As a general remark, we recall that the biological paradigm
relating L\'evy flights to the search strategies of biological organisms
is still under intense debate and the applicability of simple, attractive and ``universal''
mathematical formulations has sometimes to face complex situations
in which the details of the biological setting under consideration
may turn the problem around and require new observations, ideas and methodologies.
As a matter of fact,
the possible source of errors in the biological data is certainly multifarious. Besides the difficulty of {\em testing extensively a sufficiently large number of individuals}, possibly following
precisely all their movements in order to obtain data with a solid statistical significance, the {\em environment} in which the experiment is set may also offer, by its own nature, unexpected pitfalls for the collection of reliable results. A
 prototypical example of this difficulty is embodied by the case of wandering albatrosses who dive into the water to catch food, which, for more than a decade, in view of~\cite{NATU}, were considered as the most prominent example of forager performing L\'evy flight: in this pioneering experiment, sensors were attached to the feet of the birds and the duration of a flight was measured by the interval of time for which the sensor remained dry. In addition to the limited amount of data (five albatrosses engaged in 19 foraging bouts), as outlined in~\cite{NATUBIS, PNASS2} the interpretation of the environment may lead to different statistical interpretation of the same phenomenon (e.g. in the situation in which the albatrosses were keeping their sensor dry by sitting on an island rather than engaging a long flight).

Typically, ever-shifting environments (in which memory for the forager is of very limited use)
and situations of sparse food (in which the distance between food patches is way beyond sensory range)
may favor long-jump foraging patterns.
As a matter of fact, the influence of the environment and of the circadian cycle of the foragers' habits on different foraging patterns is also outlined in~\cite{NATU2} in the analysis of the diving depth of a blue shark. In this case, in regions with very limited diving depth, the data fit an exponential distribution, while in open ocean regions the data display power-law distributions (with an exponent~$s$ close, but\footnote{The exponent~$s=1/2$ plays often a special
role in the model since, under appropriate assumptions, it is believed to provide ``optimal''
search strategies in terms of encounter rates
corresponding to sparse, randomly distributed, nondestructive targets, see~\cite{89897765NATU}.
Though the matter has experiencing an intense debate (see~\cite{LIT1, LIT2, LIT3}), even in critical reviews the ``idea
that ``L\'evy walks can be efficient to explore space'' is overall
accepted quite broadly, together with the ``observation of power-law-like patterns in field data''
(see e.g. the last paragraph in~\cite{LIT1}). In any case, once again rigorous mathematics
can be helpful in providing unambiguous statements valid under explicit assumptions, to consolidate our knowledge of
such a difficult matter and clearly demarcate the boundaries and the limitations of our understanding.
In our paper, we do not aim at solving
the several controversies posed by the details of the L\'evy foraging hypothesis, but rather at
using one of its commonly accepted formulation to discuss the model of an ecological niche
in view of~\eqref{EQUAZIONE}--\eqref{NEUM-NLOC}. Specifically, in comparison with the controversy in~\cite{LIT1, LIT2, LIT3},
we stress that in our setting the exponent~$s$ is given and we do not optimize on it (see however~\cite{MR3590646}
for a related, but different, type of optimization problem in the fractional exponent).}
not equal, to~$1/2$). The alternation of these patterns appears to be influenced as well by the day-night cycle, since at night the shark hovered close to the surface of the sea: interestingly, if one is interested in a long-term pattern, one may try to ``average out'' the circadian cycle and fit the data by a superposition of two different distributions. 

The coexistence of L{\'e}vy and Brownian movement patterns has been tested also in~\cite{PNASS2},
also in support of the possibility that L\'evy flights may have naturally evolved in nature as an alternative beneficial search strategy ``in response to sparse resources and scant information''.

The statistical deviation from pure L\'evy flights (as well as from pure Brownian motions) has also been observed in Bumblebees: see~\cite{PhysRevLett108098103}, which ran a laboratory experiment tracking real bumblebees visiting replenishing nectar sources possibly under threat by ``artificial spiders''. When predators are present, the bumblebees perform a more careful approach to the food source to avoid spiders, hence flights with longer durations between flower visits become more frequent in presence of a predation risk.

This superposition of different statistical patterns is closely related to the notion of ``intermittent dynamics'', see e.g.~\cite{MR2299528, MR2639124, MR2670512, RevModPhys8381, Klages R}. The typical situation producing these dynamics arises from a search strategy combining phases of fast (non reactive) motion during which targets are not expected to be found, and of slow (reactive) motion in which the searcher is attentively seeking the target,
see Figure~1 in~\cite{MR2299528} (or its reproductions in Figure~1 in~\cite{RevModPhys8381} and Figure~4.5 in~\cite{Klages R}). The slow phases are typically modeled by a Brownian motion, while the fast phases take into consideration L\'evy flights (or possibly ballistic relocations in the limit as~$s\to0^+$) and often these two phases are mixed randomly.\medskip

We also recall that
a supplementary difficulty in the understanding of optimal foraging hypotheses by using {\em numerical simulations} arises from the high sensitivity of the results obtained with respect to the specific model considered (e.g., in terms of boundary conditions,
number of dimensions, uphill or downhill drifts,
see~\cite{Palyulin2931, PhysRevE78051128}).\medskip

The {\em causal interpretation} of the findings obtained is also experiencing a rather
intense debate, since it is not always evident whether
the observed L\'evy searching patterns in nature arose from an adaptive behavior or from the distributions of prey,
namely whether the anomalous diffusion in animal foraging is the outcome
of an {\em optimal search strategy} (which is itself
the byproduct of an evolutionary adaptation maximizing the
success for survival) or of an {\em interaction
between a forager and the food source distribution} -- or a combination of the two
(see e.g.~\cite{SIMS} and Section~4.3.3 in~\cite{Klages R}).
\medskip

Thus, on the one hand, the interest in the topic of animal foraging from different perspectives
bring out
the intrinsic cross-disciplinarity of the subject, which involves zoologists, ethologists, botanists, physicists, mathematicians, computer scientists, data scientists, neurologists, etc. The variety of expertise exploited is also justified
in view of the different difficulties that the topic offers (e.g., in terms of analysis of the environment, reliability of the data,
numerical implementations, neurological and evolutionary interpretations)
and the broad debate presented by a number of innovative foraging and searching strategies.
On the other hand, these
features highlight as well the importance of the mathematical deduction of coherent models from first principles (when it is possible to do so), also in order to detect the specific features which lead to different results in the analysis of the experimental data.\medskip

Specifically, we will describe here a simple model leading to~\eqref{EQUAZIONE}--\eqref{NEUM-NLOC}.
In a nutshell, these equations arise
(up to a careful tuning of the parameters)
by the combination of two dispersal occurring
with different probabilities
and by an immediate return to the niche if 
trespassing occurs.
More specifically:
\begin{itemize}
\item One can imagine that each individual is represented by
a particle subject to either a jump process (occurring with some
probability~$p$) or a random walk (occurring with some
probability~$1-p$);
\item If the particle exits the domain, it comes back immediately
to it, being redirected to any point of the domain
with the same type of probability distribution
(i.e., jumping to any point of the domain, if the egress occurred due
to the jump process, or walking to an accessible closest neighborhood,
if the egress occurred due to the random walk), and
these processes are also normalized for a unit total probability.
\end{itemize}

\begin{center}\begin{figure}
\includegraphics[width=11.9cm]{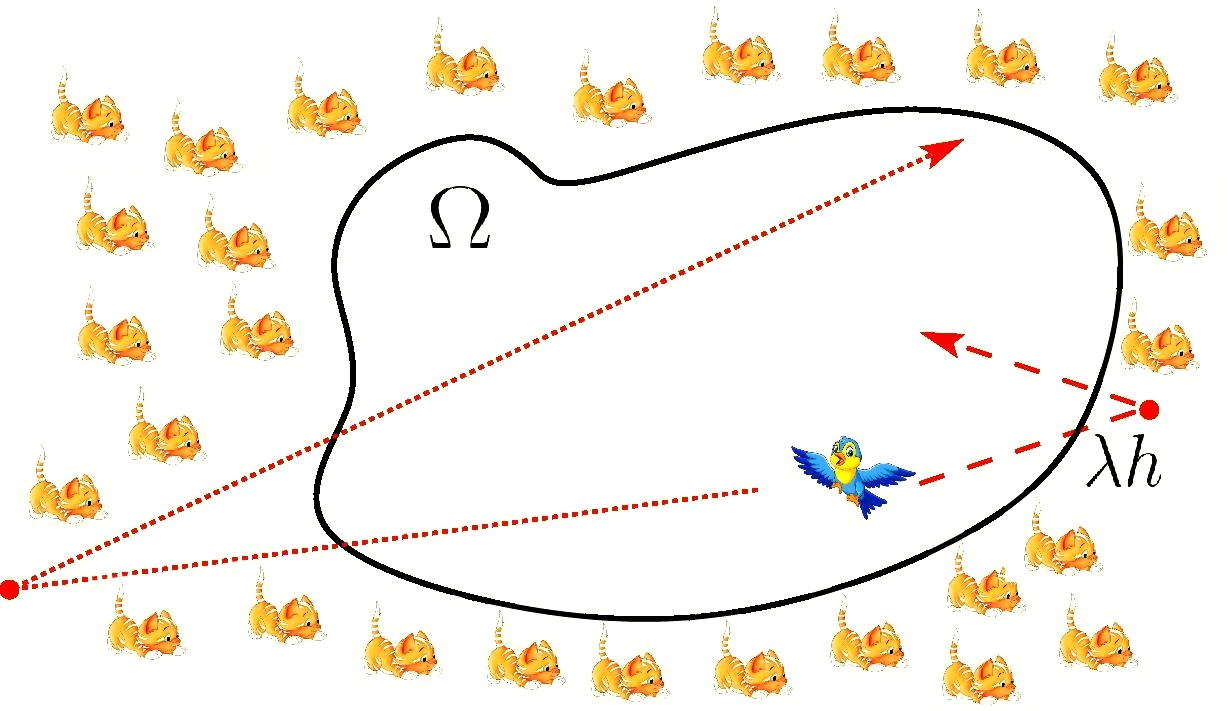}
\caption{\it Sketch of an ecological niche (the dashed
line corresponding to the classical random walk and the dotted line
to the jump process).}
\label{FIG}\end{figure}
\end{center}

See Figure~\ref{FIG} for a pictorial sketch of the notion of niche that we consider.
We observe that the link between
reflecting barriers and Neumann conditions is classical
in the theory of elliptic partial differential equations
(though the topic is often discussed only in dimension~$1$,
see~\cite[page~98]{MR2399851}).
Reflection of random walks is also classically linked
to the so-called
Skorokhod problem, see~\cite{SKORO1, SKORO2, MR837810}.
Also, the relation between jump processes
and fractional operators has been widely investigated,
see e.g.~\cite{MR2512800, MR2584076, MR3967804},
and the probabilistic methods often turn out to be extremely
profitable in the development of the analytic theory, see e.g.~\cite{MR119247, MR1438304, MR2345912}.
In our description,
for the jump process in~$\R^n$ we will exploit the
description proposed in~\cite[Section~2.1]{MR3469920},
that needs to be combined here with the idea of ``immediate
return to~$\Omega$'' that was introduced in~\cite{dirova},
and whose comprehensive probabilistic
setting has been recently provided in~\cite{VOND}.\medskip

Though built on this existing literature, the
simple model that we present here comprising both the local and the nonlocal
dispersals
is\footnote{This presentation aims at being
compelling and does not aspire to be mathematically exhaustive
(the complete proofs in a detailed probabilistic framework
having to rely on an extended version of the techniques
recently introduced in~\cite{VOND} with the goal of
comprising a classical Dirichlet form into the fractional ones).
In particular, to keep our discussion as easy 
as possible, we will often make these simplifications:
\begin{itemize}
\item
We drop the normalizing constants related to the classical
and fractional Laplacians to ease notations (this amounts to renaming
the diffusion coefficients~$\alpha$ and~$\beta$, without
affecting the arguments provided);
\item We consider here only ``formal'' expansions
(i.e., we do not estimate remainders);
\item All sets and all functions are implicitly
assumed to be nice and smooth (in particular, all functions
are supposed to be bounded and with bounded derivatives), and sequences of functions
are implicitly supposed to smoothly converge to their limits.
\end{itemize}
Typically, rigorous kernel estimates
and boundary regularity results for Neumann type
nonlocal problems rely on sophisticated analysis, see~\cite{MR4102340, AUD}.
With this, we will
obtain a
description that can be followed basically with no prerequisites,
in particular no previous knowledge of Dirichlet forms,
trace processes, Dirichlet-to-Neumann operators, or
Markov chains is needed to follow this exposition.
Also, the presentation is essentially self-contained, hence we
hope that this note can be useful for a broad range
of readers coming from different disciplines.} new.
\medskip

As a technical observation,
we remark that
we are specifically presenting here the case of indistinguishable
individuals, each presenting two diffusive patterns
with given probabilities: models describing
the case of a species with two classes of individuals
that present different diffusive patterns (such as the one
in~\cite{FALCO}) or two species with different diffusive
patterns sharing the same ecological niche (as
taken into account in~\cite{ANT}) are, in a sense, seemingly related
to our work, but they are technically quite different.

To appreciate such a difference, one can consider
the case in which local searching modes are predominant
and consider an evolution in a very large niche starting from
a very concentrated initial condition.
In this setting, models distinguishing between
individuals with short-range and long-range diffusive patterns
will find, say, at unit distance from the initial concentration point,
after a unit of time, a very small density of
short-range travelers (due to the Gaussian
decay of the classical heat equation) and a much larger density of
long-range travelers (due to the polynomial decay
of the nonlocal heat equation) and thus
the relative density of short-range travelers in that location
will be very small (even if the majority of the population
belongs to such a class). Viceversa, in the model that we consider,
in the range of parameters in which the
local search is predominant,
one can expect that short-jumps will happen with high probability
among all the individuals who have reached the location.

Interestingly, mixed diffusive patterns can arise from composite stochastic processes,
alternating intensive and area-concentrated with extensive and far-reaching search modes,
and this diffusive turnover may be induced by environmental features, such as patchy structures, see e.g.~\cite{LW}.
We also remark that diffusive operators of mixed order naturally arise in probabilistic problems
for infinitely divisible distributions via the
L\'evy-Khintchine formula and find applications in different fields of applied sciences, including phylogeny (see~\cite{FILO}) and finance (see~\cite{FINA}).
We also refer to the monograph~\cite{FORA} for
a throughout presentation
of the theory of random searches
for foraging purposes.\medskip

Moreover, from the probabilistic perspective,
we point out that the distinction between the classical (i.e., Browniann) and the jump (i.e., L\'evy)
random walks has a natural interpretation in terms of
Central Limit Theorems. Namely, the standard random walk meets the Central Limit Theorem in its classical formulation, according to which the sum of a number of independent and identically distributed random variables with finite variances converges to a normal distribution as the number of variables grows.
Instead, the jump process meets the ``generalized'' Central Limit Theorem (see~\cite{MR0233400}),
according to which the sum of a number of random variables with a power-law tails distributions having
infinite variance converges to a
stable distribution
(and this type of
stable distributions present tails that are asymptotically
proportional to the ones of a power-law distribution).
See e.g.~\cite{MR2743162} for further details and a thorough historical
review about the Central Limit Theorem,
its generalizations and its connections with stable distributions.\medskip

We also point out that anomalous diffusion driven by operators of mixed orders (possibly with
two fractional terms, as well) has been taken into consideration for the one-dimensional case without boundaries,
also in view of the Voigt profile function and of the Mellin-Barnes integrals,
and detecting two time-scales of the problem
(see equation~(3.2) in~\cite{MR2459736}
or equation~(13) in~\cite{MR2559350}, and~\cite{CRME02, CRME08, MR3886705} too).
The situation that we take into account here
takes somewhat a different point of view, since
the boundary and exterior conditions
in~\eqref{NEUM-LOC} and~\eqref{NEUM-NLOC}
play a decisive role for us in describing the mechanism of egression from, and return to, the given domain.\medskip

As a technical remark, we point out that several possible equivalent
definitions can be given for the fractional Laplacian introduced in~\eqref{EQUAZIONE}: for instance, instead of a singular integral operator, one can take into account the generator of a semigroup of operators, or Bochner's subordination formulas, or Dynkin’s formula, or Riesz potentials, or harmonic extensions, etc. -- see e.g.~\cite{MR3613319} for a comprehensive treatment of all these possible definitions.
Nevertheless, in the scientific literature, other types of fractional versions of the Laplace operators have been taken into account, in terms of censored stochastic processes and spectral analysis. These objects share some traits with the fractional Laplacian in~\eqref{EQUAZIONE} here, but they are structurally very different, see in particular~\cite{MR3233760, MR3246044} and Sections~2 and~4 in~\cite{MR4102340}.

Several reasons suggested to us that, among these alternative choices of elliptic fractional operators, the fractional Laplacian in~\eqref{EQUAZIONE} was the most convenient one to describe the random egress from a biological niche with immediate return. Indeed, the regional fractional Laplacian (see equation~(2.47) in~\cite{MR4102340}) acts on functions defined in a given domain, therefore does not seem suitable to take into account the values of the function outside the domain itself (while,
for our construction, we need to detect the exterior of the domain, where some egression takes place, though with an instantaneous
return to the domain itself).
Similarly, the spectral fractional Laplacian (see equation~(2.49) in~\cite{MR4102340}) acts on functions defined on a given domain, and more specifically on functions in the span of the eigenfunctions of the Laplacian with given boundary conditions: 
though this setting 
is adequate for diffusive problems based on the classical heat equation with suitable boundary conditions
(see~\cite{MR3796372}), not only in this case the function is not naturally defined outside the domain but also the boundary conditions are inherited directly from the local situation and do not take into account the specific features of a long-jump process
that we need to take into account in our framework.
Conversely, the fractional Laplacian in~\eqref{EQUAZIONE} has nice compatibility properties with long-jump processes
and produces an exterior condition that can be interpreted naturally as a return to the domain immediately after an egression.

We also point out that, on the one hand, the nonlocal Neumann condition in~\eqref{NEUM-NLOC}
seems to make use of a regional fractional Laplacian with respect to the domain~$\Omega$, but, on the other hand, the condition
itself takes place in the complement of~$\Omega$,  hence this prescription does not fall at once into the
conventional equations driven by the regional fractional Laplacian (though, as pointed out in~\cite{MR4102340, AUD},
a regional problem does arise by modifying the kernel with a logarithmic singularity along~$\partial\Omega$).
The delicate interplay between densities localized in the domain~$\Omega$ and quantity globally defined
in~$\R^n$ will play an important role in our construction -- for instance, the population density
is supposed to be localized in the niche, see~\eqref{DENSIINOMEPGA}, but
a convenient extension of it outside the domain will be taken into account to detect
the desired evolution equation, see~\eqref{DEFASIU0}, or equivalently one can construct 
a phantom population outside the domain by a suitable reflection method, see~\eqref{VHA}.
\medskip

In the forthcoming Section~\ref{92uiejfkn}, a more precise
description of the random process considered will be
given. {F}rom this,
we will formally derive the heat equation~\eqref{EQUAZIONE}
in Section~\ref{HEAT}
and the Neumann conditions~\eqref{NEUM-LOC}--\eqref{NEUM-NLOC}
in Section~\ref{NEUMSE}. Section~\ref{sec:concl} summarizes the conclusion of the results presented in this paper.

In the appendix, we will show an alternative approach to
the problem, in which the process is defined in the whole of~$\R^n$.
\medskip

For the reader's convenience, before diving into
the technical computations of this article,
we provide a list of the symbols and conventions
adopted
in what follows.

\section*{Notation table}

\begin{center}
\begin{tabularx}{0.8\textwidth} {
  | >{\raggedright\arraybackslash}X
  | >{\centering\arraybackslash}X
  | >{\raggedleft\arraybackslash}X | }
 \hline
{\sc Notation used }& {\sc Symbol used} & {\sc Meaning} \\
 \hline {\bf
Physical space} &$ \R^n$ & $\displaystyle{{\R\times\dots\times\R}} $ ($n$ times)\\
 \hline {\bf
Laplacian} &$ \Delta$ & $\displaystyle\sum_{j=1}^n \frac{\partial^2}{\partial x_j^2} $\\
\hline {\bf
Fractional power }& $s$ & a real number between $0$ and $1$\\
 \hline {\bf
Fractional Laplacian} &$ -(-\Delta)^s u(x) $& $\displaystyle\int_{\R^n}\frac{u(y)-u(x)}{
|y-x|^{n+2s}}\,dy $\\
 \hline {\bf
Domain } & $\Omega$  & a bounded domain in $\R^n$ with smooth boundary  \\
\hline {\bf
External unit normal }& $\nu$ & unit normal at the boundary of $\Omega$, pointing towards the exterior of $\Omega$ \\
 \hline {\bf
Lebesgue measure } & ${\mathcal{L}}^n$  & the standard $n$-dimensional measure  \\\hline{\bf
$(n-1)$-dimensional Hausdor measure}&  ${\mathcal{H}}^{n-1}$ & the standard surface measure
(for hypersurfaces of codimension $1$ in~$\R^n$) \\
 \hline {\bf
Unit ball } & $B_1$  & $\{x\in\R^n {\mbox{ s.t. }}|x|<1\}$  \\
\hline {\bf
Unit sphere } & $\partial B_1$  & $\{x\in\R^n {\mbox{ s.t. }}|x|=1\}  $\\ \hline{\bf
Indicator function of a given set $A$} & $\chi_A$  & $\chi_A(x):=\begin{cases} 1 &{\mbox{ if }}x\in A,\cr
0&{\mbox{ otherwise.}}
\end{cases} $\\
\hline {\bf
Time step }& $\tau$ & A small real number with the dimensionality of time \\
 \hline {\bf
Space step } & $h$  & A small real number with the dimensionality of space  \\
\hline {\bf
Diffusive parameters }& $\alpha$ and $\beta$ & Given positive real numbers\\
\hline {\bf
Probability of jump process to take place} & $p$ & Given real number between $0$
and $1$\\
\hline
\end{tabularx}\end{center}
\bigskip

This notation table
lists the symbols mostly used in this paper. Most of this notation is standard,
though the literature offers plenty of alternative notation
and small variances:
for instance, the Laplace operator
is sometimes denoted by~$\nabla^2$ in other papers,
the fractional Laplacian by~$|\nabla|^{2s}$
and possibly contains normalizing constants,
the normal derivative is elsewhere denoted by~$n$
(here reserved for the number of dimensions),
the infinitesimal volume element
of the Lebesgue measure by~$dx$,
the infinitesimal surface element
by~$d\sigma$, the fractional parameter~$s$ corresponds
to~$\alpha/2$ in the ``$\alpha$-stable'' literature
in probability and statistics,
and the unit sphere is elsewhere often
denoted by~${\mathbb{S}}^{n-1}$ especially in differential geometry.
We chose the notation stated in the table
since we find it more explicit and in order to circumvent unclear statements
(especially, we put an effort
in distinguishing volume and surface integrals
and in avoiding ambiguous notations related to integrals, since
the two types of measures play a significantly different role
in our framework).
Also, as a conceptual simplification,
we take our physical space to be~$\R^n$ rather than restricting
to the special cases of~$\R$, $\R^2$ and~$\R^3$
(in the case considered here,
this higher dimensional
generalization does not provide any additional difficulty
and does not make the notation heavier).

\section{Description of the random process}\label{92uiejfkn}

The probabilistic interpretation of the process leading to
equation~\eqref{EQUAZIONE}
and to the nonlocal
Neumann conditions~\eqref{NEUM-LOC}--\eqref{NEUM-NLOC}
goes as follows:
\begin{itemize}
\item $u(x, t)$ is the probability distribution of the position
of a particle moving randomly in a (say, bounded and smooth)
domain~$\Omega\subset\R^n$;
\item The
random process followed by the particle is the superposition
of a classical random walk and a long-jump
random walk. More specifically, we consider a time increment~$\tau>0$,
a space scale~$h>0$,
and an additional space parameter~$\lambda\in\N$ (in our setting,
both~$\tau$ and~$\lambda h$ will be infinitesimal);
\item
We assume that the particle
has a probability~$p$
of following a long-jump process
and a probability~$1-p$ of following a classical random walk.
The probability
for the particle of not moving at all is set to be equal to zero.
For concreteness,
we consider here the case~$p\in(0,1)$, but the cases~$p=0$
and~$p=1$ are comprised in the forthcoming discussion
(just focusing on the process allowed by such probability
and disregarding the other one);
\item When the particle exits $\Omega$,
it comes back into~$\Omega$ right away. This return to~$\Omega$
follows a natural reflection with respect to the
random process.
\end{itemize}
We now discuss the classical random walk and the
long jump process 
more precisely.

\subsection{The classical random walk}\label{CLAS1}
Concerning the classical random walk, 
we denote by~${\mathcal{L}}^{n}$ the
$n$-dimensional Lebesgue measure,
and we describe the probability of
walking from~$x\in\Omega$
to~$y\in\Omega$ as the superposition of
these two phenomena:
\begin{itemize}
\item the particle can walk ``directly''
from~$x\in\Omega$ to~$y\in\Omega$,
which occurs with
probability density~$\frac{\chi_{[0,\lambda h)}(|x-y|)}{{\mathcal{L}}^{n}(B_{\lambda h})}$:
that is, each points in the ball~$B_{\lambda h}(x)$ are reachable from~$x$
with uniform distribution;
\item the particle can also walk
from~$x$ to~$y$ ``after a reflection''
in the complement of~$\Omega$: namely,
the particle can walk from~$x$ to~$
z\in\R^n\setminus\Omega$, 
which
occurs with probability
density~$\frac{\chi_{[0,\lambda h)}(|x-z|)}{{\mathcal{L}}^{n}(B_{\lambda h})}$, and then
walk back instantaneously from~$
z$ to~$y\in\Omega$,
and this occurs with probability
density~$
\frac{\chi_{[0,\lambda h)}(|y-z|)}{{\mathcal{L}}^{n}(\Omega\cap
B_{\lambda h}(z))}$.
That is, the random walk can reach points outside the domain,
from which it bounces back to the domain with the same
uniform probability law.
As a result, the probability
density of a reflected walk
from~$x\in\Omega$
to~$y\in\Omega$ is equal to
$$ 
\int_{\R^n\setminus\Omega}
\frac{\chi_{[0,\lambda h)}(|x-z|)
\;\chi_{[0,\lambda h)}(|y-z|)}{
{\mathcal{L}}^{n}(B_{\lambda h})\;{\mathcal{L}}^{n}(\Omega\cap
B_{\lambda h}(z))}\,dz.
$$
\end{itemize}
The total probability density 
of a walk
from~$x\in\Omega$
to~$y\in\Omega$ is therefore the sum of the probabilities of a direct walk
and a reflected walk
and it is equal to
\begin{equation}\label{PW} P_W(x\to y):=
\frac{\chi_{[0,\lambda h)}(|x-y|)}{{\mathcal{L}}^{n}(B_{\lambda h})}+
\int_{\R^n\setminus\Omega}
\frac{\chi_{[0,\lambda h)}(|x-z|)
\;\chi_{[0,\lambda h)}(|y-z|)}{ {\mathcal{L}}^{n}(B_{\lambda h})\;
{\mathcal{L}}^{n}(\Omega\cap
B_{\lambda h}(z))}\,dz
.\end{equation}
We remark that the probability of walking back to the domain after
exiting is equal to~$1$, since, for every~$z\in\R^n\setminus\Omega$,
$$
\int_{\Omega}
\frac{\chi_{[0,\lambda h)}(|y-z|)}{ 
{\mathcal{L}}^{n}(\Omega\cap
B_{\lambda h}(z))}\,dy
=1.$$
Also, the probability of walking from a given point
of~$\Omega$ to the domain~$\Omega$ itself is equal to~$1$,
since, for each~$x\in\Omega$,
\begin{equation}\label{COIS}
\begin{split}&
\int_{\Omega}
P_W(x\to y)\,dy\\=\,&
\int_\Omega\frac{\chi_{[0,\lambda h)}(|x-y|)}{{\mathcal{L}}^{n}(B_{\lambda h})}\,dy
+\int_\Omega\left[
\int_{\R^n\setminus\Omega}
\frac{\chi_{[0,\lambda h)}(|x-z|)
\;\chi_{[0,\lambda h)}(|y-z|)}{ {\mathcal{L}}^{n}(B_{\lambda h})\;
{\mathcal{L}}^{n}(\Omega\cap
B_{\lambda h}(z))}\,dz\right]\,dy\\=\,&
\frac{ {\mathcal{L}}^{n}(B_{\lambda h}(x)\cap\Omega)}{{\mathcal{L}}^{n}(B_{\lambda h})}
+\int_{\R^n\setminus\Omega}\left[
\int_\Omega
\frac{\chi_{[0,\lambda h)}(|x-z|)
\;\chi_{[0,\lambda h)}(|y-z|)}{ {\mathcal{L}}^{n}(B_{\lambda h})\;
{\mathcal{L}}^{n}(\Omega\cap
B_{\lambda h}(z))}\,dy\right]\,dz\\=\,&
\frac{ {\mathcal{L}}^{n}(B_{\lambda h}(x)\cap\Omega)}{{\mathcal{L}}^{n}(B_{\lambda h})}
+\int_{\R^n\setminus\Omega}
\frac{\chi_{[0,\lambda h)}(|x-z|)}{ {\mathcal{L}}^{n}(B_{\lambda h})} \,dz\\=\,&
\frac{ {\mathcal{L}}^{n}(B_{\lambda h}(x)\cap\Omega)}{
{\mathcal{L}}^{n}(B_{\lambda h})}
+
\frac{ {\mathcal{L}}^{n}(B_{\lambda h}(x)\setminus\Omega)}{
{\mathcal{L}}^{n}(B_{\lambda h})}\\=\,&1.
\end{split}\end{equation}
We also remark that
$$ P_W(y\to x)=P_W(x\to y),$$
thanks to~\eqref{PW}.
This and~\eqref{COIS} entail that,
for every~$\eta\in\Omega$,
\begin{equation}\label{EPPS}
\int_{ \Omega}
P_W(\xi\to \eta)\,d\xi
=
\int_{ \Omega}
P_W(\eta\to \xi)\,d\xi=1.
\end{equation}

\subsection{The jump process}\label{CLAS1J}

Concerning the jump process, 
we denote by~${\mathcal{H}}^{n-1}$
the $(n-1)$-dimensional Hausdorff measure
and
we let, for all~$z\in\R^n\setminus\Omega$,
\begin{equation}\label{2.4BIS}
\mu_h(z):=2s \,h^{2s} \int_{\Omega\setminus B_h(z)}
\frac{dy}{|y-z|^{n+2s}} .
\end{equation}
Given~$x$, $y\in\Omega$, 
the probability
of jumping from~$x$
to~$y$ consists in the superposition of two events:
\begin{itemize}
\item the particle can jump ``directly''
from~$x$
to~$y$: 
this occurs with probability density
$$\frac{2s \,h^{2s}\,\chi_{(h,+\infty)}(|x-y|)}{{\mathcal{H}}^{n-1}(\partial B_1)\;|x-y|^{n+2s}} ;$$
\item the particle can jump
from~$x$
to~$y$ ``after a reflection'' in the complement of~$\Omega$:
that is, the particle can jump from~$x$ to any point~$
z\in\R^n\setminus\Omega$,
which occurs with probability
density
$$\frac{2s \,h^{2s}\,\chi_{(h,+\infty)}(|x-z|)}{{\mathcal{H}}^{n-1}(\partial B_1)\;|x-z|^{n+2s}} ,$$
and then
instantaneously back to~$\Omega$, that is,
from~$z$
to~$y\in\Omega$,
which occurs with probability density
$$\frac{2s \,h^{2s}\,\chi_{(h,+\infty)}(|y-z|)}{
\mu_h(z)\;|y-z|^{n+2s}} .$$
Consequently,
the probability density of a reflected jump
from~$x$
to~$y$ passing through~$z\in\R^n\setminus\Omega$ is
obtained by the product rule of independent events, and it is equal to
$$\int_{\R^n\setminus\Omega}\frac{(2s)^2 \,h^{4s}\,\chi_{(h,+\infty)}(|x-z|)
\,\chi_{(h,+\infty)}(|y-z|)
}{{\mathcal{H}}^{n-1}(\partial B_1)\;\mu_h(z)\;|x-z|^{n+2s}|y-z|^{n+2s}} \,dz.$$
\end{itemize}
The probability density
of a jump
from~$x$
to~$y$ is therefore the sum of the probability densities of a direct jump
and a reflected jump
and it is equal to
\begin{equation}\label{7687697654397532} \begin{split}&
P_J(x\to y):=
\frac{2s \,h^{2s}\,\chi_{(h,+\infty)}(|x-y|)}{{\mathcal{H}}^{n-1}(\partial B_1)\;|x-y|^{n+2s}}
+\int_{\R^n\setminus\Omega}\frac{(2s)^2 \,h^{4s}\,\chi_{(h,+\infty)}(|x-z|)
\,\chi_{(h,+\infty)}(|y-z|)
}{{\mathcal{H}}^{n-1}(\partial B_1)\;\mu_h(z)\;|x-z|^{n+2s}|y-z|^{n+2s}} \,dz.\end{split}
\end{equation}
We remark that the probability of jumping back to the domain after
exiting is equal to~$1$, since, for each~$
z\in\R^n\setminus\Omega$,
$$ \int_\Omega
\frac{2s \,h^{2s}\,\chi_{(h,+\infty)}(|y-z|)}{
\mu_h(z)\;|y-z|^{n+2s}}\,dy=1,$$
thanks to~\eqref{2.4BIS}. In this sense, the return to the domain
follows the same jump law of the egress, with the term~$\mu_h(z)$
acting as a normalization probability factor.

Also, the probability of jumping from a given point
of~$\Omega$ to the domain~$\Omega$ itself is equal to~$1$,
since, for each~$x\in\Omega$,
\begin{eqnarray*}&&
\int_{\Omega}P_J(x\to y)\,dy\\&=&
\int_{\Omega}
\frac{2s \,h^{2s}\,\chi_{(h,+\infty)}(|x-y|)}{{\mathcal{H}}^{n-1}(\partial B_1)\;|x-y|^{n+2s}}
\,dy+\int_{\Omega}\left[
\int_{\R^n\setminus\Omega}\frac{(2s)^2 \,h^{4s}\,\chi_{(h,+\infty)}(|x-z|)
\,\chi_{(h,+\infty)}(|y-z|)
}{{\mathcal{H}}^{n-1}(\partial B_1)\;\mu_h(z)\;|x-z|^{n+2s}|y-z|^{n+2s}} \,dz
\right]\,dy
\\&=&\frac{2s \,h^{2s}}{{\mathcal{H}}^{n-1}(\partial B_1)}
\int_{\Omega\setminus B_h(x)}
\frac{dy}{|x-y|^{n+2s}}
+\int_{\R^n\setminus\Omega}\left[
\int_{\Omega}\frac{(2s)^2 \,h^{4s}\,\chi_{(h,+\infty)}(|x-z|)
\,\chi_{(h,+\infty)}(|y-z|)
}{{\mathcal{H}}^{n-1}(\partial B_1)\;\mu_h(z)\;|x-z|^{n+2s}|y-z|^{n+2s}} \,dy
\right]\,dz
\\&=&\frac{2s \,h^{2s}}{{\mathcal{H}}^{n-1}(\partial B_1)}
\int_{\Omega\setminus B_h(x)}
\frac{dy}{|x-y|^{n+2s}}
+\int_{\R^n\setminus\Omega}
\frac{2s \,h^{2s}\,\chi_{(h,+\infty)}(|x-z|)
}{{\mathcal{H}}^{n-1}(\partial B_1)\;|x-z|^{n+2s}}\,dz
\\&=&\frac{2s \,h^{2s}}{{\mathcal{H}}^{n-1}(\partial B_1)}\left[
\int_{\Omega\setminus B_h(x)}
\frac{dy}{|x-y|^{n+2s}}
+\int_{(\R^n\setminus\Omega)\setminus B_h(x)}
\frac{dz
}{|x-z|^{n+2s}} \right]\\
\\&=&\frac{2s \,h^{2s}}{{\mathcal{H}}^{n-1}(\partial B_1)}
\int_{\R^n \setminus B_h(x)}
\frac{d\xi}{|x-\xi|^{n+2s}}
\\&=& 2s \,h^{2s}
\int_{h}^{+\infty}
\frac{d\rho}{\rho^{1+2s}}\\&=&1.
\end{eqnarray*}
Since~$P_J(y\to x)=P_J(x\to y)$, due to~\eqref{7687697654397532},
we also have that, for all~$x\in\Omega$,
\begin{equation}\label{FolAMS}
\int_{\zeta\in\Omega}P_J(\zeta\to x)\,d{\mathcal{H}}^{n-1}_\zeta=1.
\end{equation}
Following~\cite{VERO}, we consider here
the special sets of parameters
\begin{equation}\label{CA-01}
\tau:=h^{2s},\qquad\lambda:=h^{s-1}\in\N, \qquad{\mbox{
and~$p$ is independent of $h$}}
\end{equation}
and we describe in detail how these
prescriptions lead to a heat equation of mixed local and nonlocal
type as in~\eqref{EQUAZIONE}, with nonlocal Neumann conditions
as in~\eqref{NEUM-LOC}--\eqref{NEUM-NLOC}.
These aspects will be discussed in the forthcoming sections.

\section{The heat equation}\label{HEAT}

We denote by~$u(x,t)$
the total probability density for the particle to be at the point~$x$
at time~$t$, and we set
$$ u_J(x,t):=p \,u(x,t)\qquad{\mbox{ and }}\qquad
u_W(x,t):=(1-p) \,u(x,t).$$
Roughly speaking, one can consider~$u_W$ (respectively, $u_J$) to be
the contribution of probability density for the particle to be at the point~$x\in \Omega$
at time~$t\in\tau\N$ coming from the classical random walk
(respectively, the jump process).
By construction, both the classical random walk and
the jump process do not leave the domain~$\Omega$, hence we can write that
\begin{equation}\label{DENSIINOMEPGA}
u(x,t)=u_W(x,t)=u_J(x,t)=0\qquad{\mbox{
whenever }}\,x\not\in\Omega.
\end{equation}
Moreover,
\begin{equation}\label{7yhbCA50-1} u(x,t)=u_J(x,t)+ u_W(x,t).\end{equation}
Given~$x\in\Omega$,
the probability of being at~$x$ at time~$t+\tau$ is
the superposition of the probability of being at any other point~$y\in\Omega$ at time~$t$, times the probability~$P(y\mapsto x)$ of jumping 
from~$y$ to~$x$. For this reason, we write
\begin{equation}\label{7yhbCA50-2}
u(x,t+\tau)=\int_{\Omega} u(y,t)\,P(y\mapsto x)\,dy.\end{equation}
By the random processes described in Sections~\ref{CLAS1}
and~\ref{CLAS1J}, we have that
\begin{equation}\label{7yhbCA50-3} P(y\mapsto x)=p \,P_J(y\mapsto x)+(1-p)\,P_W(y\mapsto x),\end{equation}
and consequently,
by~\eqref{EPPS} and~\eqref{FolAMS},
$$ \int_{\Omega } P(y\mapsto x)\,dy=1.$$
This and~\eqref{7yhbCA50-2} lead to
\begin{equation*} \frac{u(x,t+\tau)-u(x,t)}\tau
=\int_{\Omega} \frac{u(y,t)-u(x,t)}\tau\,P(y\mapsto x)\,dy.\end{equation*}
As a result,
by~\eqref{7yhbCA50-1}
and~\eqref{7yhbCA50-3},
\begin{equation}\label{8u76u9i6030f984}
\begin{split}&
\left[ 
\frac{u_W(x,t+\tau)-u_W(x,t)}\tau-
(1-p)\int_{\Omega} \frac{u(y,t)-u(x,t)}\tau\,P_W(y\mapsto x)\,dy\right]\\
&\qquad+
\left[ 
\frac{u_J(x,t+\tau)-u_J(x,t)}\tau-
p\int_{\Omega} \frac{u(y,t)-u(x,t)}\tau\,P_J(y\mapsto x)\,dy\right]
\\=\,&
\frac{u(x,t+\tau)-u(x,t)}\tau-
\int_{\Omega} \frac{u(y,t)-u(x,t)}\tau\,P(y\mapsto x)\,dy=0.
\end{split}
\end{equation}
We denote by~${\mathcal{T}}_W(x,t)$ (respectively, ${\mathcal{T}}_J(x,t)$)
the item in the first (respectively, the second) square brackets
in~\eqref{8u76u9i6030f984}.

Recalling~\eqref{PW}, we know that,
for every~$x\in\Omega$ with\footnote{As customary,
we use the notation~``$\Subset$'' to mean ``compactly contained'',
that is, given two subsets~$A$ and~$B$ of~$\R^n$, we say that~$A\Subset B$ if~$A$
is bounded and the closure of~$A$ is contained in~$B$ (equivalently, if there exists a compact set~$K$
such that~$A\subseteq K\subseteq B$).}~$B_{\lambda h}(x)\Subset\Omega$,
\begin{eqnarray*}
P_W(x\to y)&=&
\frac{\chi_{[0,\lambda h)}(|x-y|)}{{\mathcal{L}}^{n}(B_{\lambda h})}+
\int_{\R^n\setminus\Omega}
\frac{\chi_{[0,\lambda h)}(|x-z|)
\;\chi_{[0,\lambda h)}(|y-z|)}{ {\mathcal{L}}^{n}(B_{\lambda h})\;
{\mathcal{L}}^{n}(\Omega\cap
B_{\lambda h}(z))}\,dz
\\&=&\frac{\chi_{[0,\lambda h)}(|x-y|)}{{\mathcal{L}}^{n}(B_{\lambda h})},
\end{eqnarray*}
since the integration set over~$z$ is empty.

Thus,
for every~$x\in\Omega$ with~$B_{\lambda h}(x)\Subset\Omega$,
\begin{equation}\label{7uJHNSNDu93uytg}
{\mathcal{T}}_W(x,t)=
\frac{u_W(x,t+\tau)-u_W(x,t)}{\tau}-
\frac{1-p}{{\mathcal{L}}^{n}( B_{\lambda h})}\int_{B_{\lambda h}(x)}\frac{
u(y,t)-u(x,t)}{\tau}\,dy.
\end{equation}
Furthermore, using that~$\lambda h=\sqrt\tau$, due to~\eqref{CA-01},
if~$y\in B_{\lambda h}(x)$, a formal Taylor expansion leads to
\begin{eqnarray*} u(y,t)&=&u(x,t)+\nabla u(x,t)\cdot(y-x)
+\frac12\,D^2u(x,t)(x-y)\cdot(x-y)+O(|x-y|^3)\\&=&
u(x,t)+\nabla u(x,t)\cdot(y-x)+\frac12\,\sum_{i,j=1}^n
\partial^2_{ij}u(x,t)(x_i-y_i)(x_j-y_j)
+O(\tau^{3/2}).\end{eqnarray*}
This, together with two odd symmetry cancellations, gives that
\begin{eqnarray*}&&
\int_{B_{\lambda h}(x)}\frac{
u(y,t)-u(x,t)}{\tau}\,dy\\&=&\frac1\tau
\int_{B_{\sqrt\tau}(x)}\left[
\nabla u(x,t)\cdot(y-x)+\frac12\,\sum_{i,j=1}^n
\partial^2_{ij}u(x,t)(x-y)\cdot(x-y)\right]
\,dy+O\left(\tau^{\frac{n+1}2}\right)\\&=&\sum_{i,j=1}^n\frac{\partial^2_{ij}u(x,t)}{2\tau}
\int_{B_{\sqrt\tau}(x)}
(x_i-y_i)(x_j-y_j)
\,dy+O\left(\tau^{\frac{n+1}2}\right)\\&=&\sum_{i,j=1}^n\frac{\tau^{\frac{n}2}\partial^2_{ij}u(x,t)}{2}
\int_{B_1}
\omega_i\omega_j
\,d\omega+O\left(\tau^{\frac{n+1}2}\right)\\&=&\sum_{i=1}^n\frac{\tau^{\frac{n}2}\partial^2_{i}u(x,t)}{2}
\int_{B_1}
\omega_i^2
\,d\omega+O\left(\tau^{\frac{n+1}2}\right).
\end{eqnarray*}
Consequently, since, for each~$i\in\{1,\dots,n\}$,
$$
c_o:=\int_{B_1}
|\omega|^2\,d\omega
=
\sum_{j=1}^n\int_{B_1}
\omega_j^2\,d\omega=
n\int_{B_1}
\omega_i^2
\,d\omega,$$
we find that
$$ \int_{B_{\lambda h}(x)}\frac{
u(y,t)-u(x,t)}{\tau}\,dy=
\sum_{i=1}^n\frac{c_o\;\tau^{\frac{n}2}\partial^2_{i}u(x,t)}{2n}
+O\left(\tau^{\frac{n+1}2}\right)
= \frac{c_o\;\tau^{\frac{n}2}\Delta u(x,t)}{2n}
+O\left(\tau^{\frac{n+1}2}\right)
.$$
Thus, plugging this information into~\eqref{7uJHNSNDu93uytg},
we obtain
\begin{eqnarray*}
{\mathcal{T}}_W(x,t)&=&
\frac{u_W(x,t+\tau)-u_W(x,t)}{\tau}-
\frac{1-p}{{\mathcal{L}}^{n}( B_{\sqrt\tau})}
\left[
\frac{c_o\;\tau^{\frac{n}2}\Delta u(x,t)}{2n}
+O\left(\tau^{\frac{n+1}2}\right)
\right]
\\&=&\frac{u_W(x,t+\tau)-u_W(x,t)}{\tau}-
\frac{1-p}{{\mathcal{L}}^{n}( B_{1})}
\left[
\frac{c_o\;\Delta u(x,t)}{2n}
+O(\sqrt\tau)
\right].
\end{eqnarray*}
Passing to the limit the previous identity, we thereby formally
conclude that, for all~$x\in\Omega$,
\begin{equation}\label{SDKMdr0ot4}
\lim_{h\to0}
{\mathcal{T}}_W(x,t)=\partial_t u_W(x,t)-
(1-p)\,\Delta u(x,t),
\end{equation}
up to neglecting normalization constants.

As for the contributions coming from the jump process,
recalling~\eqref{CA-01}, we see that
\begin{equation}\label{vf34rfB5thnUY06YSh}
{\mathcal{T}}_J(x,t)=
\frac{u_J(x,t+\tau)-u_J(x,t)}\tau-
p\int_{\Omega} \frac{u(y,t)-u(x,t)}{h^{2s}}\,P_J(y\mapsto x)\,dy.\end{equation}
Moreover, in light of~\eqref{7687697654397532},
\begin{equation}\label{prv0oov4}
\begin{split}&\int_{\Omega} \frac{u(y,t)-u(x,t)}{h^{2s}}\,P_J(y\mapsto x)\,dy\\
=\;&
\int_{\Omega} \frac{u(y,t)-u (x,t)}{h^{2s}}\,
\Bigg\{
\frac{2s \,h^{2s}\,\chi_{(h,+\infty)}(|x-y|)}{{\mathcal{H}}^{n-1}(\partial B_1)\;|x-y|^{n+2s}}
\\&\qquad+\int_{\R^n\setminus\Omega}\frac{(2s)^2 \,h^{4s}\,\chi_{(h,+\infty)}(|x-z|)
\,\chi_{(h,+\infty)}(|y-z|)
}{{\mathcal{H}}^{n-1}(\partial B_1)\;\mu_h(z)\;|x-z|^{n+2s}|y-z|^{n+2s}} \,dz\Bigg\}
\,dy\\=\;&
\frac{2s}{{\mathcal{H}}^{n-1}(\partial B_1)}\,\big(I_h+J_h\big),
\end{split}\end{equation}
where
\begin{eqnarray*}
I_h&:=&
\int_{\Omega\setminus B_h(x)} \frac{u(y,t)-u (x,t)}{|x-y|^{n+2s}}\,
dy
\\{\mbox{and }}\quad J_h&:=&2s
\,h^{2s}\int_\Omega\left[
\int_{\R^n\setminus\Omega}\frac{\big(u(y,t)-u (x,t)\big)\;\chi_{(h,+\infty)}(|x-z|)
\,\chi_{(h,+\infty)}(|y-z|)
}{\mu_h(z)\;|x-z|^{n+2s}|y-z|^{n+2s}} \,dz\right]\,dy
.\end{eqnarray*}
It is now convenient to set
\begin{equation}\label{DEFASIU0} \begin{split}
&U(x,t):=\begin{dcases}
u(x,t) & {\mbox{ if }}x\in\overline\Omega,\\
\\
\frac{\displaystyle\int_\Omega \frac{u (y,t)}{|x-y|^{n+2s}}\,dy}{
\displaystyle\int_\Omega \frac{dy}{|x-y|^{n+2s}}
}& {\mbox{ if }}x\in\R^n\setminus\overline\Omega\end{dcases}\end{split}\end{equation}
and
\begin{equation}\label{Titagscrane} \begin{split}
&U_h(x,t):=\begin{dcases}
u(x,t) & {\mbox{ if }}x\in\overline\Omega,\\
\\
\frac{\displaystyle\int_{\Omega \setminus B_h(x)}\frac{u (y,t)}{|x-y|^{n+2s}}\,dy}{
\displaystyle\int_{\Omega \setminus B_h(x)}\frac{dy}{|x-y|^{n+2s}}
}& {\mbox{ if }}x\in\R^n\setminus\overline\Omega.\end{dcases}\end{split}
\end{equation}
In some sense, the setting
in~\eqref{DEFASIU0}
can be seen
as a ``nonlocal variant'' of the classical ``method of images''
(see e.g. page~28 in~\cite{MR3380662} or Chapter~2 in~\cite{ELEME}), in which
one defines a ``phantom'' population
outside the niche, to keep the population
balance constant in the niche.

With this notation, and recalling the definition of~$\mu_h$ in~\eqref{2.4BIS},
we observe that, for all~$x\in\Omega$,
\begin{eqnarray*}
J_h&:=&
\int_{(\R^n\setminus\Omega)\setminus B_h(x)}\left[
\int_{\Omega\setminus B_h(z)}\frac{ u(y,t)-u (x,t)
}{\displaystyle\int_{\Omega\setminus B_h(z)}\frac{dw}{|w-z|^{n+2s}}\;
|y-z|^{n+2s}} \,dy\right]\,\frac{
dz}{|x-z|^{n+2s}}
\\&=&
\int_{(\R^n\setminus\Omega)\setminus B_h(x)}\big(U_h(z,t)-U_h(x,t)\big)\,\frac{
dz}{|x-z|^{n+2s}}.
\end{eqnarray*}
Hence, since
\begin{equation}
\label{FN-con}
U_{ h }(x,t)\longrightarrow U(x,t)\qquad\qquad{\mbox{as }}\,h\to0,
\end{equation}
we formally obtain that
$$ \lim_{h\to0}
(I_h+J_h)=\int_{\R^n}\big(U(z,t)-U(x,t)\big)\,\frac{
dz}{|x-z|^{n+2s}}=-(-\Delta)^s U(x,t),$$
up to neglecting normalization constants.

{F}rom this and~\eqref{prv0oov4}, we deduce that
$$ \lim_{h\to0}
\int_{\Omega} \frac{u(y,t)-u(x,t)}{h^{2s}}\,P_J(y\mapsto x)\,dy
=-(-\Delta)^s U(x,t),$$
up to neglecting normalization constants once again.

This and~\eqref{vf34rfB5thnUY06YSh} lead to
\begin{equation*}
\lim_{h\to0}
{\mathcal{T}}_J(x,t)=\partial_tu_J(x,t)+p(-\Delta)^s U (x,t).
\end{equation*}
{F}rom this, \eqref{8u76u9i6030f984} and~\eqref{SDKMdr0ot4},
and omitting constants again, we find, when~$x\in\Omega$,
\begin{eqnarray*}0&=&
\lim_{h\to0}
{\mathcal{T}}_W(x,t)+{\mathcal{T}}_J(x,t)
\\&=&
\partial_t u_W(x,t)-(1-p)
\Delta u(x,t)+\partial_tu_J(x,t)+p(-\Delta)^s U (x,t)\\
&=&\partial_t u(x,t)-(1-p)
\Delta u(x,t)+p(-\Delta)^s U (x,t)
\\&=&\partial_t U(x,t)-(1-p)
\Delta U(x,t)+p(-\Delta)^s U (x,t),
\end{eqnarray*}
which is the desired heat equation in~\eqref{EQUAZIONE}.

\section{The Neumann condition}\label{NEUMSE}

We suppose, up to a translation, that the origin lies in an $h$-neighborhood
of the boundary of~$\Omega$ -- in fact,
for simplicity, let us just suppose that~$0\in\partial\Omega$
and let~$\nu$ be the external unit normal of~$\Omega$.

The probability density of finding the particle at~$x\in\Omega$ at time~$t+\tau$
can be written as in~\eqref{7yhbCA50-2}, thus leading to~\eqref{8u76u9i6030f984}.
We will therefore exploit~\eqref{8u76u9i6030f984} in this setting, that we write in the form
\begin{equation}\label{8u76u9i6030f984-INSA}
\begin{split}&
\left[ 
\frac{u_W(x,t+\tau)-u_W(x,t)}{\lambda h}-
(1-p)\int_{\Omega} \frac{u(y,t)-u(x,t)}{\lambda h}\,P_W(y\mapsto x)\,dy\right]\\
&\qquad+
\left[ 
\frac{u_J(x,t+\tau)-u_J(x,t)}{\lambda h}-
p\int_{\Omega} \frac{u(y,t)-u(x,t)}{\lambda h}\,P_J(y\mapsto x)\,dy\right]=0.
\end{split}
\end{equation}
We denote by~${\mathcal{S}}_W(x,t)$ (respectively, ${\mathcal{S}}_J(x,t)$)
the item in the first (respectively, the second)
square brackets
in~\eqref{8u76u9i6030f984-INSA}.

Our goal is now to exploit suitable
expansions of the function~$u(x,t)$
to relate it
to the mixed fractional equation~\eqref{EQUAZIONE}.
For this, we will relate this function to stable densities,
and specifically the term~${\mathcal{S}}_W$
will correspond to a local diffusive operator
while the term~${\mathcal{S}}_J$ to a fractional one.
Indeed, the Brownian probability term~$P_W$,
acting on a short range and being invariant under rotations
and translations, will produce in the limit the Laplace operator,
which constitutes the 
Gaussian contribution of the evolution equation in~\eqref{EQUAZIONE}.
Instead, the power-law probability term~$P_J$
will maintain its long tail in the limit and produce the fractional
Laplace operator,
which constitutes the 
L\'evy contribution of~\eqref{EQUAZIONE}.

More precisely, we notice
that, if~$y\in \overline{B_{2\lambda h}(x)}$,
with a formal Taylor expansion we can write that
$$ u(y,t)-u(x,t)=\nabla u(x,t)\cdot(y-x)+O(|x-y|^2)=\nabla u(x,t)\cdot(y-x)+O(\lambda^2 h^2).$$
We also remark that, due to the construction of the random walk in Section~\ref{CLAS1},
the function~$y\mapsto P_W(y\mapsto x)$ is supported in~$y\in \overline{B_{2\lambda h}(x)}$, and, as a result,
for all~$y\in\R^n$ we have that
$$ \big(u(y,t)-u(x,t)\big)\,P_W(y\mapsto x)=\Big(\nabla u(x,t)\cdot(y-x)+O(\lambda^2 h^2)\Big)
P_W(y\mapsto x).$$
This gives that
$$ {\mathcal{S}}_W(x,t)=
\frac{u_W(x,t+\tau)-u_W(x,t)}{\lambda h}-
(1-p)\int_{\Omega} \left(\frac{\nabla u(x,t)\cdot(y-x)}{\lambda h}+O(\lambda h)\right)\,P_W(y\mapsto x)\,dy.
$$
Also, in light of~\eqref{CA-01}, another formal Taylor expansion gives that
\begin{equation}\label{GEb9n6hbcoqac} \frac{u_W(x,t+\tau)-u_W(x,t)}{\lambda h}=\frac{O(\tau)}{\lambda h}
=O(h^{s}),\end{equation}
and therefore
\begin{equation}\label{ESSEFASew} {\mathcal{S}}_W(x,t)=
-
(1-p)\int_{\Omega} \frac{\nabla u(x,t)\cdot(y-x)}{\lambda h}\,P_W(y\mapsto x)\,dy+O(h^s).\end{equation}
As a consequence, recalling~\eqref{PW},
\begin{equation}\label{23Jc}\begin{split}
{\mathcal{S}}_W(x,t)=\,&
-\frac{1-p}{{\mathcal{L}}^{n}(B_{1})}
\int_{\Omega} \frac{\nabla u(x,t)\cdot(y-x)}{(\lambda h)^{n+1}}
\,\Bigg[\chi_{[0,\lambda h)}(|x-y|)\\&\qquad\qquad+
\int_{\R^n\setminus\Omega}
\frac{\chi_{[0,\lambda h)}(|x-z|)
\;\chi_{[0,\lambda h)}(|y-z|)}{ 
{\mathcal{L}}^{n}(\Omega\cap
B_{\lambda h}(z))}\,dz
\Bigg]\,dy+O(h^s)\\=\,&-\frac{1-p}{{\mathcal{L}}^{n}(B_{1})}\Bigg[
\int_{\Omega\cap B_{\lambda h}(x)} \frac{\nabla u(x,t)\cdot(y-x)}{
(\lambda h)^{n+1}}\,dy
\\&\qquad+
\int_{(\R^n\setminus\Omega)\cap B_{\lambda h}(x)}
\left(\int_{\Omega\cap B_{\lambda h}(z)} \frac{\nabla u(x,t)\cdot(y-x)}{
(\lambda h)^{n+1}\,{\mathcal{L}}^{n}(\Omega\cap
B_{\lambda h}(z))}\,dy\right)\,dz
\Bigg]+O(h^s)\\=\,&-\frac{1-p}{{\mathcal{L}}^{n}(B_{1})}\Bigg[
\int_{\frac{\Omega-x}{\lambda h}\cap B_1} \nabla u(x,t)\cdot Y\,dY
\\&\qquad+
\int_{\frac{(\R^n\setminus\Omega)-x}{\lambda h}\cap B_1}
\left(\int_{\Big(\frac{\Omega-x}{\lambda h}-Z\Big)\cap B_1} \frac{\nabla u(x,t)\cdot (Y+Z)}{
{\mathcal{L}}^{n}\left(\Big({\frac{\Omega-x}{\lambda h}-Z\Big)\cap B_1}\right)}\,dY\right)\,dZ
\Bigg]+O(h^s)
.
\end{split}\end{equation}
Now we consider the jump process and we exploit~\eqref{GEb9n6hbcoqac}
(with~$u_J$ in place of~$u_W$), thus obtaining that
\begin{equation}\label{7say3e9rjv9retisyay345}
{\mathcal{S}}_J(x,t)=-
p\int_{\Omega} \frac{u(y,t)-u(x,t)}{\lambda h}\,P_J(y\mapsto x)\,dy+O(h^s).
\end{equation}
Recalling the definition of~$\mu_h$ in~\eqref{2.4BIS}
and that of~$U_h$ in~\eqref{Titagscrane},
we also observe that
\begin{eqnarray*}E_h&:=&
\int_{\Omega}\left[
\int_{\R^n\setminus\Omega}\frac{2s\,h^{4s}\,\chi_{(h,+\infty)}(|x-z|)
\,\chi_{(h,+\infty)}(|y-z|)\;\big(u(y,t)-u(x,t)\big)
}{\lambda h\;\mu_h(z)\;|x-z|^{n+2s}|y-z|^{n+2s}} \,dz\right]\,dy\\&=&
\int_{(\R^n\setminus\Omega)\setminus B_h(x)}\left[\int_{\Omega\setminus B_h(z)}
\frac{h^{2s}\,\big(u(y,t)-u(x,t)\big)
}{\lambda h\;\displaystyle\int_{\Omega\setminus B_h(z)}
\frac{dw}{|w-z|^{n+2s}}\;|x-z|^{n+2s}|y-z|^{n+2s}} \,dy\right]\,dz\\&=&
\int_{(\R^n\setminus\Omega)\setminus B_h(x)}
\frac{h^{2s}\,\big(U_h(z,t)-U_h(x,t)\big)
}{\lambda h\;|x-z|^{n+2s}} \,dz.
\end{eqnarray*}
Formally assuming sufficient regularity for~$U_h$,
we write that
$$ \big|U_h(z,t)-U_h(x,t)\big|\le C\,\min\{1,|x-z|\}$$
for some~$C>0$, whence, by the parameters' choice in~\eqref{CA-01},
\begin{eqnarray*}
|E_h|&\le& C\,h^s\int_{\R^n\setminus B_h(x)}
\frac{\min\{1,|x-z|\}
}{|x-z|^{n+2s}} \,dz\\&=&
C\,h^s\int_{B_1(x)\setminus B_h(x)}\frac{
dz
}{|x-z|^{n+2s-1}} +O(h^s)\\&=&O(h^{1-s})+O(h^s)+O(h^{1/2}|\ln h|).
\end{eqnarray*}
{F}rom this and~\eqref{7687697654397532},
and possibly renaming~$C>0$,
we deduce that
\begin{eqnarray*}&&
\left|\int_{\Omega} \frac{u(y,t)-u(x,t)}{\lambda h}\,P_J(y\mapsto x)\,dy\right|\le
\int_{\Omega\setminus B_h(x)}
\frac{2s \,h^{2s}\,\big|u(y,t)-u(x,t)\big|}{\lambda h\;{\mathcal{H}}^{n-1}(\partial B_1)\;|x-y|^{n+2s}}
\,dy
+\frac{2s\;|E_h|}{{\mathcal{H}}^{n-1}(\partial B_1)}\\&&\qquad\le
Ch^s\int_{\Omega\setminus B_h(x)}
\frac{dy}{|x-y|^{n+2s-1}}
+O(h^{1-s})+O(h^s)+O(h^{1/2}|\ln h|)\\&&\qquad
=O(h^{1-s})+O(h^s)+O(h^{1/2}|\ln h|).
\end{eqnarray*}
As a consequence, we infer from \eqref{7say3e9rjv9retisyay345} that
\[ {\mathcal{S}}_J(x,t)=O(h^{1-s})+O(h^{s})+O(h^{1/2}|\ln h|).\]
{F}rom this, \eqref{8u76u9i6030f984-INSA} and~\eqref{23Jc}
we get that
\begin{equation*}
\begin{split}&
\int_{\frac{\Omega-x}{\lambda h}\cap B_1} \nabla U(x,t)\cdot Y\,dY
+\int_{\frac{(\R^n\setminus\Omega)-x}{\lambda h}\cap B_1}
\left(\int_{\Big(\frac{\Omega-x}{\lambda h}-Z\Big)\cap B_1} \frac{\nabla U(x,t)\cdot (Y+Z)}{
{\mathcal{L}}^{n}\left(\Big({\frac{\Omega-x}{\lambda h}-Z\Big)\cap B_1}\right)}
\,dY\right)\,dZ
\\&\qquad\qquad\qquad=
O(h^{1-s})+O(h^{s})+O(h^{1/2}|\ln h|).
\end{split}\end{equation*}
In particular, taking the formal limit as~$x\to0$,
\begin{equation}\label{0kdf-565}
\begin{split}&
\int_{\frac{\Omega}{\lambda h}\cap B_1} \nabla U(0,t)\cdot Y\,dY
+\int_{\frac{\R^n\setminus\Omega}{\lambda h}\cap B_1}
\left(\int_{\Big(\frac{\Omega}{\lambda h}-Z\Big)\cap B_1} \frac{\nabla U(0,t)\cdot (Y+Z)}{
{\mathcal{L}}^{n}\left(\Big({\frac{\Omega}{\lambda h}-Z\Big)\cap B_1}\right)}\,dY\right)\,dZ
\\&\qquad\qquad\qquad=
O(h^{1-s})+O(h^{s})+O(h^{1/2}|\ln h|).
\end{split}\end{equation}
Now we observe that, as~$h\to0$, the rescaled
domain~$\frac{\Omega}{\lambda h}$
approaches the
halfspace~$\Pi$ through the origin with external normal~$\nu(0)$.

As a result, we take the formal limit as~$h\to0$ of the identity
in~\eqref{0kdf-565}, gathering that
\begin{equation}\label{09km-02erfev0}
\int_{\Pi\cap B_1} \nabla U(0,t)\cdot Y\,dY
+\int_{(\R^n\setminus\Pi)\cap B_1}
\left(\int_{\Big(\Pi-Z\Big)\cap B_1} \frac{\nabla U(0,t)\cdot (Y+Z)}{
{\mathcal{L}}^{n}\left(\Big({\Pi-Z\Big)\cap B_1}\right)}\,dY\right)\,dZ=0
.\end{equation}
Now we claim that
\begin{equation}\label{09km-02erfev}
\int_{\Pi\cap B_1} Y\,dY
+\int_{(\R^n\setminus\Pi)\cap B_1}
\left(\int_{\Big(\Pi-Z\Big)\cap B_1} \frac{Y+Z}{
{\mathcal{L}}^{n}\left(\Big({\Pi-Z\Big)\cap B_1}\right)}\,dY\right)\,dZ=-
c_\star\,\nu(0),
\end{equation}
for some~$c_\star>0$.
To prove this,
up to a rotation, we suppose that~$\nu(0)=(0,\dots,0,1)$,
whence~$\Pi=\{ \zeta=(\zeta_1,\dots,\zeta_n)\in\R^n{\mbox{ s.t. }}\zeta_n<0\}$.
Therefore, the $n$th component of the left hand side of~\eqref{09km-02erfev}
is
\begin{eqnarray*}
&&\int_{\{Y_n<0\}\cap B_1} Y_n\,dY
+\int_{\{Z_n>0\}\cap B_1}
\left(\int_{\Big(\{Y_n<0\}-Z\Big)\cap B_1} \frac{Y_n+Z_n}{
{\mathcal{L}}^{n}\left(\Big({\{Y_n<0\}-Z\Big)\cap B_1}\right)}\,dY\right)\,dZ,
\end{eqnarray*}
which is strictly negative, say equal to~$-c_\star$ for some~$c_\star>0$,
since both the numerators in the above integrals are negative.

Also, for each~$i\in\{1,\dots,n-1\}$, the~$i$th component of the left hand side of~\eqref{09km-02erfev}
is
\begin{eqnarray*}
&&\int_{\{Y_n<0\}\cap B_1} Y_i\,dY
+\int_{\{Z_n>0\}\cap B_1}
\left(\int_{\Big(\{Y_n<0\}-Z\Big)\cap B_1} \frac{Y_i+Z_i}{
{\mathcal{L}}^{n}\left(\Big({\{Y_n<0\}-Z\Big)\cap B_1}\right)}\,dY\right)\,dZ,
\end{eqnarray*}
which is equal to zero, since the integrands are odd with respect
to the $i$th variable, while the domains remain the same after switching~$Y_i$
to~$-Y_i$ and~$Z_i$ to~$-Z_i$.

These considerations lead to
\begin{eqnarray*}&&
\int_{\Pi\cap B_1} Y\,dY
+\int_{(\R^n\setminus\Pi)\cap B_1}
\left(\int_{\Big(\Pi-Z\Big)\cap B_1} \frac{Y+Z}{
{\mathcal{L}}^{n}\left(\Big({\Pi-Z\Big)\cap B_1}\right)}\,dY\right)\,dZ\\&&\qquad
=(0,\dots,0,-c_\star)=-c_\star\,\nu(0),
\end{eqnarray*}
and this proves~\eqref{09km-02erfev}, as desired.

By~\eqref{09km-02erfev0}
and~\eqref{09km-02erfev} we deduce that
$$ -c_\star\,\nabla U(0,t)\cdot\nu(0)=0,$$
that is~$\partial_\nu U(0,t)=0$. Since the origin was an arbitrary point of~$\partial\Omega$,
we have thus obtained the local Neumann prescription in~\eqref{NEUM-LOC}.

Furthermore, 
using~\eqref{DEFASIU0},
we see that
the nonlocal Neumann prescription in~\eqref{NEUM-NLOC}
is fulfilled in~$\R^n\setminus\overline\Omega$,
as desired.

\section{Conclusions}\label{sec:concl}

We have proposed
a mathematical description of an ecological model, by also providing
concrete motivations.
The model has been described precisely at an intuitive level,
with explicit calculations that can be followed without
many prerequisites.

The model describes 
a biological population which can pursue both a
classical (i.e. Brownian) and a long-range (i.e. L\'evy)
dispersal strategies and 
lives in an ecological niche. In this situation, 
the combination of Brownian and L\'evy processes
produces an evolution equation 
similar to the heat flow but with the Laplace operator
replaced by the superposition of a Laplacian and a fractional
Laplacian.

The equation is set in a domain (which corresponds
to the biological niche).
Additionally, in view of this combination of different types of diffusion,
a new type of Neumann condition is needed to impose
a zero-flux prescription of the niche and model the return
to the niche after a possible egress due to a stochastic process.
This new type of Neumann condition is also a superposition
of the classical Neumann condition (prescribing the normal
derivative at the boundary of the domain) and of a fractional
Neumann condition (corresponding to an integral relation
in the complement of the domain).

The analytic approach to this model can provide inspiration
for further understanding of biological situations of interest.

\begin{appendix}
\section{The ``point of view of the niche''}

It is interesting to point out that the function~$U$
introduced in~\eqref{DEFASIU0} gives a possible
different interpretation of the stochastic process,
providing an equivalent setting defined in the whole of~$\R^n$
and not only in~$\Omega$ (remarkably, also in~\cite{VOND}
a rigorous approach is provided to define a
random process in~$\R^n$
equivalent to the one producing the Neumann
condition of~\cite{dirova}).

That is, one could define
\begin{equation}\label{DP-2}d\pi_W(y;x):=
\frac{\chi_{[0,\lambda h)}(|x-y|)}{{\mathcal{L}}^{n}(B_{\lambda h})}\,dy
\end{equation}
and
\begin{equation}\label{DP} d\pi_J(y;x):=
\frac{2s \,h^{2s}\,\chi_{(h,+\infty)}(|x-y|)}{{\mathcal{H}}^{n-1}
(\partial B_1)\;|x-y|^{n+2s}}\,dy.
\end{equation}
In a sense, the quantity in~\eqref{DP-2}
is induced
by the random walk probability 
density described in 
Section~\ref{CLAS1}, and
the quantity in~\eqref{DP} is coming
from the jump process probability 
density described in 
Section~\ref{CLAS1J}.
It is also worth pointing out that
the processes in~\eqref{DP-2} and~\eqref{DP} somewhat
``complement'' each other, since
the first takes care of short-range movements at distance
less than~$\lambda h$ and the latter produces
trajectories that are always longer than~$h$.

We will
now complement these objects by
a suitable ``method of images''.
Namely, we set
\begin{equation}\label{TOPSO} d\pi(y;x):=p\,d\pi_J(y;x)+
(1-p)\,d\pi_W(y;x)\end{equation}
and
$$ \Omega^{(\lambda h)}:=
\{ x\in \Omega{\mbox{ s.t. }}
B_{\lambda h}(x)\Subset\Omega\},$$
and we could consider~$V$
as the density of a biological
population in~$ \Omega$
and that of a ``phantom'' population~$V_h$ outside~$ \Omega^{(\lambda h)}$,
assuming that such a population is subject
to a random
process, with probability density of jump from~$y\in\R^n$
to~$x\in\Omega^{(\lambda h)}$
given by~\eqref{TOPSO};
the phantom population is then
reset to the appropriate value in~\eqref{DEFASIU0}
outside~$\Omega^{(\lambda h)}$ at each time step of the process.
Namely, while the process occurs
in~$\Omega^{(\lambda h)}$, the diffusing density is defined outside~$\Omega^{(\lambda h)}$ by
\begin{equation}\label{VHA} V_h(x,t):=\begin{dcases}
V(x,t) &{\mbox{ if }}x\in\Omega^{(\lambda h)},\\
\\
\displaystyle\frac{\displaystyle
\int_{\Omega\cap B_{\lambda h}(x)} V(y,t)\,dy}{
\displaystyle{\mathcal{L}}^n(
\Omega\cap B_{\lambda h}(x))}
&{\mbox{ if }}x\in\overline\Omega\setminus\Omega^{(\lambda h)},\\
\\
\frac{\displaystyle\int_\Omega \frac{V (y,t)}{|x-y|^{n+2s}}\,dy}{
\displaystyle\int_\Omega \frac{dy}{|x-y|^{n+2s}}}
&{\mbox{ if }}x\in\R^n\setminus\overline\Omega.
\end{dcases}\end{equation}
Suggestively, this definition is
somehow a minor variation of the setting in~\eqref{DEFASIU0},
and we will make the formal ansatz that the external population~$V_h$
smoothly converges as~$h\to0$, calling, with a slight abuse of notation,
$V$ also this limit function defined outside~$\Omega$.

We remark that, differently from Section~\ref{CLAS1J},
the process here takes place in~$\R^n$, since also
the ``phantom'' population is subject to it, being forced
to jump in~$\Omega^{(\lambda h)}$.
Interestingly, this ``phantom'' population is
chosen precisely to fit with both the
local and nonlocal Neumann prescriptions.
In this spirit, for each~$x\in\Omega^{(\lambda h)}$,
the random process leads to
\begin{equation}\label{NOPOe4S0}
V_h(x,t+\tau)=\int_{\R^n} V_h(y,t)\,d\pi(y;x).
\end{equation}
It is interesting to observe that
\begin{equation*}
\begin{split}\int_{\R^n}
d\pi_J(y;x)=\int_{\R^n}
\frac{2s \,h^{2s}\,\chi_{(h,+\infty)}(|x-y|)}{
{\mathcal{H}}^{n-1}(\partial B_1)\;|x-y|^{n+2s}}=
2s \,h^{2s}\,\int_h^{+\infty}\frac{d\rho}{\rho^{1+2s}}=1
\end{split}\end{equation*}
and
\begin{equation*}
\int_{\R^n}d\pi_W(y;x)=\int_{\R^n}
\frac{\chi_{[0,\lambda h)}(|x-y|)}{{\mathcal{L}}^{n}(B_{\lambda h})}\,dy
=\frac{{\mathcal{L}}^{n}(B_{\lambda h}(x))}{{\mathcal{L}}^{n}(B_{\lambda h})}
=1,\end{equation*}
which lead to
\begin{equation}\label{NOPOe4S}
\int_{\R^n}
d\pi(y;x)=1.\end{equation}
It is suggestive to observe that~\eqref{NOPOe4S}
does not really say that~$d\pi(y;x)$ is
a probability density, since it is measuring
the jumps from~$y\in\R^n$ to~$x\in\Omega^{(\lambda h)}$ and not viceversa.
As a matter of fact,
the ``phantom'' population
balance is set
as in~\eqref{VHA} to keep the population
constant in the niche.
{F}rom a different perspective,
we may consider~\eqref{NOPOe4S}
the ``point of view of the niche'', which, at any unit of time,
is expected to receive a bit of biological population
(either the ``real'' population coming from the niche
or the ``phantom'' one).

By~\eqref{NOPOe4S0} and~\eqref{NOPOe4S},
it follows that
\begin{eqnarray*}&&
V_h(x,t+\tau)-V_h(x,t)=
\int_{\R^n} \big(V_h(y,t)-V_h(x,t)\big)\,d\pi(y;x)\\&&\qquad=
\frac{2sp \,h^{2s}}{
{\mathcal{H}}^{n-1}(\partial B_1)}
\int_{\R^n\setminus B_h(y)}
\frac{\big(V_h(y,t)-V_h(x,t)\big)}{|x-y|^{n+2s}}\,dy+\frac{1-p}{{\mathcal{L}}^{n}(B_{\lambda h})}
\int_{B_{\lambda h}(x)}
\big(V_h(y,t)-V_h(x,t)\big)\,dy.
\end{eqnarray*}
Hence, recalling the parameter choice in~\eqref{CA-01},
dividing by~$
\tau=h^{2s}=(\lambda h)^2$, and then sending~$h\to0$, we formally obtain
\begin{eqnarray*}&&
\partial_tV(x,t)\\&=&\lim_{h\to0}
\Bigg( \frac{2sp}{
{\mathcal{H}}^{n-1}(\partial B_1)}
\int_{\R^n\setminus B_h(y)}
\frac{\big(V_h(y,t)-V_h(x,t)\big)}{|x-y|^{n+2s}}\,dy\\&&\qquad\qquad\qquad
+\frac{1-p}{(\lambda h)^2{\mathcal{L}}^{n}(B_{\lambda h})}
\int_{B_{\lambda h}(x)}
\big(V_h(y,t)-V_h(x,t)\big)\,dy
\Bigg)\\&=&-p(-\Delta)^sV(x,t)\\&&\quad
+\lim_{h\to0}\frac{1-p}{(\lambda h)^2{\mathcal{L}}^{n}(B_{\lambda h})}
\int_{B_{\lambda h}(x)}
\left(\nabla V_h(x,t)\cdot(y-x)+\frac12\,D^2V_h(x,t)(y-x)\cdot(y-x)\right)\,dy
\\&=&-p(-\Delta)^sV(x,t)\\&&\quad
+\lim_{h\to0}\frac{1-p}{2{\mathcal{L}}^{n}(B_{1})}
\int_{B_{1}}D^2V_h(x,t) z\cdot z\,dz\\&=&-p(-\Delta)^sV(x,t)+(1-p)\Delta V(x,t),
\end{eqnarray*}
for all~$x\in\Omega$, up to normalization constants that we omit,
and this corresponds to the diffusive equation in~\eqref{SDKMdr0ot4}.

Furthermore, the setting in~\eqref{VHA}
formal recovers the classical and nonlocal
Neumann prescriptions
in~\eqref{NEUM-LOC}
and~\eqref{NEUM-NLOC}.
Indeed, assuming that~$V_h\to V$ as~$h\to0$, we deduce from~\eqref{VHA}
that, for all~$x\in\R^n\setminus\overline\Omega$,
$$ V(x,t)=\frac{\displaystyle\int_\Omega \frac{V (y,t)}{|x-y|^{n+2s}}\,dy}{
\displaystyle\int_\Omega \frac{dy}{|x-y|^{n+2s}}},$$
which recovers the nonlocal Neumann condition in~\eqref{NEUM-NLOC}.

As for the classical Neumann condition, we take~$p\in\partial\Omega$
and check this condition at~$p$. Up to a translation, we can suppose that~$p$
is the origin. Then, since~$0\in\partial\Omega\subset\overline\Omega$,
we deduce from~\eqref{VHA} that
\begin{equation}\label{VH1} V_h(0,t)=\frac{\displaystyle
\int_{\Omega\cap B_{\lambda h}} V(y,t)\,dy}{
\displaystyle{\mathcal{L}}^n(
\Omega\cap B_{\lambda h})}.\end{equation}
Also, if~$\nu(0)$ is the exterior normal of~$\Omega$ at the
origin and~$M\ge2$, we have that~$x_{h,M}:=-\frac{\nu(0)\,\lambda h}M\in 
\Omega\setminus\Omega^{(\lambda h)}$ if~$h>0$ is small enough,
and thus~\eqref{VHA} gives that
\begin{equation}\label{VH2} V_h(x_{h,M},t)=\frac{\displaystyle
\int_{\Omega\cap B_{\lambda h}(x_{h,M})} V(y,t)\,dy}{
\displaystyle{\mathcal{L}}^n(
\Omega\cap B_{\lambda h}(x_{h,M}))}.\end{equation}
Also, taking a formal limit,
\begin{eqnarray*}&&
\lim_{h\to0}\frac{V_h(x_{h,M},t)-V_h(0,t)}{\lambda h}=
\lim_{h\to0}\frac{\nabla V_h(0,t)\cdot x_{h,M}+O(\lambda ^2 h^2)}{\lambda h}=\lim_{h\to0}
\left( -\frac{\nu(0)}M\cdot\nabla V_h(0,t)+O(\lambda h)\right)\\&&\qquad=
-\frac{\nu(0)}M\cdot\nabla V(0,t)=-\frac{\partial_\nu V(0,t)}M.
\end{eqnarray*}
Combining this with~\eqref{VH1} and~\eqref{VH2}, we write that
\begin{align}\label{SPbspmne94}&
-\frac{\partial_\nu V(0,t)}M\\=\,&\lim_{h\to0}\frac1{\lambda h}\left(
\frac{\displaystyle\int_{\Omega\cap B_{\lambda h}(x_{h,M})} V(y,t)\,dy}{
\displaystyle{\mathcal{L}}^n(\Omega\cap B_{\lambda h}(x_{h,M}))}
-\frac{\displaystyle\int_{\Omega\cap B_{\lambda h}} V(y,t)\,dy}{
\displaystyle{\mathcal{L}}^n(
\Omega\cap B_{\lambda h})}
\right)\\ =\,&
\lim_{h\to0}\frac1{\lambda h}\left(
\frac{\displaystyle\int_{\Omega\cap B_{\lambda h}(x_{h,M})} \Big(V(y,t)-V(0,t)\Big)\,dy}{
\displaystyle{\mathcal{L}}^n(\Omega\cap B_{\lambda h}(x_{h,M}))}
-\frac{\displaystyle\int_{\Omega\cap B_{\lambda h}} \Big(V(y,t)-V(0,t)\Big)\,dy}{
\displaystyle{\mathcal{L}}^n(
\Omega\cap B_{\lambda h})}
\right)\\ =\,&
\lim_{h\to0}\frac1{\lambda h}\left(
\frac{\displaystyle\int_{\Omega\cap B_{\lambda h}(x_{h,M})} \nabla V(0,t)\cdot y\,dy+O(\lambda^{n+2}h^{n+2})}{
\displaystyle{\mathcal{L}}^n(\Omega\cap B_{\lambda h}(x_{h,M}))}
-\frac{\displaystyle\int_{\Omega\cap B_{\lambda h}} \nabla V(0,t)\cdot y\,dy+O(\lambda^{n+2}h^{n+2})}{
\displaystyle{\mathcal{L}}^n(
\Omega\cap B_{\lambda h})}
\right)\\ =\,&
\lim_{h\to0}\frac1{\lambda h}\left(
\frac{\displaystyle\int_{\Omega\cap B_{\lambda h}(x_{h,M})} \nabla V(0,t)\cdot y\,dy}{
\displaystyle{\mathcal{L}}^n(\Omega\cap B_{\lambda h}(x_{h,M}))}
-\frac{\displaystyle\int_{\Omega\cap B_{\lambda h}} \nabla V(0,t)\cdot y\,dy}{
\displaystyle{\mathcal{L}}^n(
\Omega\cap B_{\lambda h})}
\right)\\ =\,&
\lim_{h\to0}\left(
\frac{\displaystyle\int_{\frac{\Omega-x_{h,M}}{\lambda h}
\cap B_{1}} \nabla V(0,t)\cdot \left(z-\frac{\nu(0)}M\right)\,dz}{
\displaystyle{\mathcal{L}}^n\left(\frac{\Omega-x_{h,M}}{\lambda h}
\cap B_{1}\right)}
-\frac{\displaystyle\int_{\frac{\Omega}{\lambda h}\cap B_{1}} \nabla V(0,t)\cdot z\,dz}{
\displaystyle{\mathcal{L}}^n\left(\frac{\Omega}{\lambda h}\cap B_{1}\right)}
\right)
.
\end{align}
Up to a rotation, we can assume that the exterior normal of~$\Omega$ at the
origin is~$(0,\dots,0,1)$, hence~$\frac\Omega{\lambda h}$
approaches the halfspace~$\Pi:=\{z=(z_1,\dots,z_n)\in\R^n {\mbox{ s.t. }}z_n<0\}$
as~$h\to0$. Therefore,~$\frac{\Omega-x_{h,M}}{\lambda h}=
\frac{\Omega}{\lambda h}+\frac{\nu(0)}M$ approaches~$\Pi+\frac{(0,\dots,0,1)}{M}$.
{F}rom these observations and~\eqref{SPbspmne94},
we obtain that
\begin{equation}\label{76878675670987}
-\frac{\partial_\nu V(0,t)}M
=
\frac{\displaystyle\int_{\left(\{z_n<0\}
+\frac{(0,\dots,0,1)}{M}\right)
\cap B_{1}} \nabla V(0,t)\cdot \left(z-\frac{(0,\dots,0,1)}M\right)\,dz}{
\displaystyle{\mathcal{L}}^n\left(\left(\{z_n<0\}
+\frac{(0,\dots,0,1)}{M}\right)
\cap B_{1}\right)}
-\frac{\displaystyle\int_{\{z_n<0\}\cap B_{1}} \nabla V(0,t)\cdot z\,dz}{
\displaystyle{\mathcal{L}}^n\left(\{z_n<0\}\cap B_{1}\right)}
.\end{equation}
We also remark that
$$ \{z_n<0\}
+\frac{(0,\dots,0,1)}{M}=\left\{z_n<\frac1M\right\}$$
and thus
$$ {\mathcal{L}}^n \left(\left(\{z_n<0\}
+\frac{(0,\dots,0,1)}{M}\right)\cap B_1\right)=
a_0+\frac{\varpi}{M}+o\left(\frac{1}{M}\right),$$
where~$a_0:={\mathcal{L}}^n(\{z_n<0\}\cap B_1)$
and~$\varpi$ is the $(n-1)$-dimensional
Lebesgue measure of the unit ball in~$\R^{n-1}$.
As a result,
$$ \frac{1}{{\mathcal{L}}^n \left(\left(\{z_n<0\}
+\displaystyle\frac{(0,\dots,0,1)}{M}\right)\cap B_1\right)}=
\frac{1}{a_0\left(1+\displaystyle\frac{\varpi}{a_0M}+o\left(\displaystyle\frac{1}{M}\right)\right)}=
\frac{1-\displaystyle\frac{\varpi}{a_0M}+o\left(\displaystyle\frac{1}{M}\right)}{a_0}
.$$
This and~\eqref{76878675670987} give that
\begin{equation*}
\begin{split}&
-\frac{a_0\,\partial_\nu V(0,t)}M\\
=\,&
\displaystyle\int_{\left\{z_n<\frac{1}{M}\right\}
\cap B_{1}} \nabla V(0,t)\cdot \left(z-\frac{(0,\dots,0,1)}M\right)\,dz
\left[1-\displaystyle\frac{\varpi}{a_0M}+o\left(\displaystyle\frac{1}{M}\right)\right]\\&\quad
-\displaystyle\int_{\{z_n<0\}\cap B_{1}} \nabla V(0,t)\cdot z\,dz\\
=\,&
\displaystyle\int_{\left\{z_n<\frac{1}{M}\right\}
\cap B_{1}} \nabla V(0,t)\cdot z\,dz
\left[1-\displaystyle\frac{\varpi}{a_0M}\right]-\frac1M
\displaystyle\int_{\left\{z_n<\frac{1}{M}\right\}
\cap B_{1}} \partial_\nu V(0,t)\,dz
\\&\quad
-\displaystyle\int_{\{z_n<0\}\cap B_{1}} \nabla V(0,t)\cdot z\,dz
+o\left(\displaystyle\frac{1}{M}\right)\\=\,&
\displaystyle\int_{\left\{0<z_n<\frac{1}{M}\right\}\cap B_{1}} \nabla V(0,t)\cdot z\,dz
-\displaystyle\frac{\varpi}{a_0M}
\displaystyle\int_{\left\{z_n<\frac{1}{M}\right\}
\cap B_{1}} \nabla V(0,t)\cdot z\,dz\\&\quad
-\frac{a_0\, \partial_\nu V(0,t)}M
+o\left(\displaystyle\frac{1}{M}\right).
\end{split}\end{equation*}
This and a symmetric cancellation in the first $n-1$
variables give that
\begin{eqnarray*}
0&=&
\displaystyle\int_{\left\{0<z_n<\frac{1}{M}\right\}\cap B_{1}} \partial_n V(0,t)\, z_n\,dz
-\displaystyle\frac{\varpi}{a_0M}
\displaystyle\int_{\left\{z_n<0\right\}
\cap B_{1}} \partial_n V(0,t)\, z_n\,dz
+o\left(\displaystyle\frac{1}{M}\right)\\&=&
-\displaystyle\frac{\varpi\,\partial_\nu V(0,t)}{a_0M}
\displaystyle\int_{\left\{z_n<0\right\}
\cap B_{1}} z_n\,dz
+o\left(\displaystyle\frac{1}{M}\right).
\end{eqnarray*}
Hence, letting
$$ b_0:=-\int_{\left\{z_n<0\right\}
\cap B_{1}} z_n\,dz,$$
we notice that~$b_0>0$ and we conclude that
$$ 0=\displaystyle\frac{\varpi\,b_0
\partial_\nu V(0,t)}{a_0M}+o\left(\displaystyle\frac{1}{M}\right).$$
Multiplying by~$M$ and sending~$M\to+\infty$, we thereby find that~$
\varpi\,b_0
\partial_\nu V(0,t)=0$, and consequently~$\partial_\nu V(0,t)=0$.
This is precisely the classical Neumann condition,
and the proposed strategy to obtain it can be seen
as a further ``nonlocal variant'' of the classical ``method of images''
(see e.g. page~28 in~\cite{MR3380662}).
\end{appendix}


\begin{bibdiv}
\begin{biblist}

\bib{MR4102340}{article}{
   author={Abatangelo, Nicola},
   title={A remark on nonlocal Neumann conditions for the fractional
   Laplacian},
   journal={Arch. Math. (Basel)},
   volume={114},
   date={2020},
   number={6},
   pages={699--708},
   issn={0003-889X},
   doi={10.1007/s00013-020-01440-9},
}

\bib{MR3967804}{article}{
   author={Abatangelo, Nicola},
   author={Valdinoci, Enrico},
   title={Getting acquainted with the fractional Laplacian},
   conference={
      title={Contemporary research in elliptic PDEs and related topics},
   },
   book={
      series={Springer INdAM Ser.},
      volume={33},
      publisher={Springer, Cham},
   },
   date={2019},
   pages={1--105},
}

\bib{MR2512800}{book}{
   author={Applebaum, David},
   title={L\'{e}vy processes and stochastic calculus},
   series={Cambridge Studies in Advanced Mathematics},
   volume={116},
   edition={2},
   publisher={Cambridge University Press, Cambridge},
   date={2009},
   pages={xxx+460},
   isbn={978-0-521-73865-1},
   doi={10.1017/CBO9780511809781},
}

\bib{AUD}{article}{
       author = {Audrito, Alessandro},
       author = {Felipe-Navarro, Juan-Carlos},
       author = {Ros-Oton, Xavier},
        title = {The Neumann problem for the fractional Laplacian: regularity up to the boundary},
      journal = {arXiv e-prints},
date = {2020},
          eid = {arXiv:2006.10026},
        pages = {arXiv:2006.10026},
archivePrefix = {arXiv},
      adsnote = {Provided by the SAO/NASA Astrophysics Data System}
}

\bib{MR2345912}{article}{
   author={Ba\~{n}uelos, Rodrigo},
   author={Bogdan, Krzysztof},
   title={L\'{e}vy processes and Fourier multipliers},
   journal={J. Funct. Anal.},
   volume={250},
   date={2007},
   number={1},
   pages={197--213},
   issn={0022-1236},
   doi={10.1016/j.jfa.2007.05.013},
}

\bib{LW}{article}{
ISSN={0012-9658},
volume={88},
issue={8},
date={2007},
title={How many animals really do the L\'evy walk?},
pages={1962–1969},
author={Benhamou, Simon},
journal={Ecology},
doi={10.1890/06-1769.1},
}

\bib{MR2299528}{article}{
   author={B\'{e}nichou, O.},
   author={Coppey, M.},
   author={Moreau, M.},
   author={Voituriez, R.},
   title={Intermittent search strategies: when losing time becomes
   efficient},
   journal={Europhys. Lett.},
   volume={75},
   date={2006},
   number={2},
   pages={349--354},
   issn={0295-5075},
   doi={10.1209/epl/i2006-10100-3},
}

\bib{RevModPhys8381}{article}{
  title = {Intermittent search strategies},
  author = {B\'enichou, O.},
  author = {Loverdo, C.},
  author = {Moreau, M.}
  author = {Voituriez, R.},
  journal = {Rev. Mod. Phys.},
  volume = {83},
  issue = {1},
  pages = {81--129},
  numpages = {0},
  year = {2011},
  doi = {10.1103/RevModPhys.83.81},
  url = {https://link.aps.org/doi/10.1103/RevModPhys.83.81}
}

\bib{FINA}{article}{
   author={Bingham, N. H.},
   author={Kiesel, R\"{u}diger},
   title={Semi-parametric modelling in finance: theoretical foundations},
   journal={Quant. Finance},
   volume={2},
   date={2002},
   number={4},
   pages={241--250},
   issn={1469-7688},
   doi={10.1088/1469-7688/2/4/201},
}

\bib{MR119247}{article}{
   author={Blumenthal, R. M.},
   author={Getoor, R. K.},
   title={Some theorems on stable processes},
   journal={Trans. Amer. Math. Soc.},
   volume={95},
   date={1960},
   pages={263--273},
   issn={0002-9947},
   doi={10.2307/1993291},
}

\bib{MR1438304}{article}{
   author={Bogdan, Krzysztof},
   title={The boundary Harnack principle for the fractional Laplacian},
   journal={Studia Math.},
   volume={123},
   date={1997},
   number={1},
   pages={43--80},
   issn={0039-3223},
   doi={10.4064/sm-123-1-43-80},
}

\bib{MR3469920}{book}{
   author={Bucur, Claudia},
   author={Valdinoci, Enrico},
   title={Nonlocal diffusion and applications},
   series={Lecture Notes of the Unione Matematica Italiana},
   volume={20},
   publisher={Springer, [Cham]; Unione Matematica Italiana, Bologna},
   date={2016},
   pages={xii+155},
   isbn={978-3-319-28738-6},
   isbn={978-3-319-28739-3},
   doi={10.1007/978-3-319-28739-3},
}

\bib{LIT2}{article}{
  title = {Comment on ``Inverse Square L\'evy Walks are not Optimal Search Strategies for $d{\ge}2$''},
  author = {Buldyrev, S. V.},
    author = {Raposo, E. P.},
      author = {Bartumeus, F.},
        author = {Havlin, S.},
          author = {Rusch, F. R.},
            author = {da Luz, M. G. E.},
              author = {Viswanathan, G. M.},
  journal = {Phys. Rev. Lett.},
  volume = {126},
  issue = {4},
  pages = {048901},
  numpages = {2},
  year = {2021},
  doi = {10.1103/PhysRevLett.126.048901},
  url = {https://link.aps.org/doi/10.1103/PhysRevLett.126.048901}
}

\bib{CRME08}{article}{
  title = {Generalized fractional diffusion equations for accelerating subdiffusion and truncated L\'evy flights},
  author = {Chechkin, A. V.},
    author = {Gonchar, V. Yu.},
      author = {Gorenflo, R.},
        author = {Korabel, N.},
          author = {Sokolov, I. M.},
  journal = {Phys. Rev. E},
  volume = {78},
  issue = {2},
  pages = {021111},
  numpages = {13},
  year = {2008},
  doi = {10.1103/PhysRevE.78.021111},
  url = {https://link.aps.org/doi/10.1103/PhysRevE.78.021111}
}

\bib{CRME02}{article}{
  title = {Retarding subdiffusion and accelerating superdiffusion governed by distributed-order fractional diffusion equations},
  author = {Chechkin, A. V.},
    author = {Gorenflo, R.},
      author = {Sokolov, I. M.},
  journal = {Phys. Rev. E},
  volume = {66},
  issue = {4},
  pages = {046129},
  numpages = {7},
  year = {2002},
  doi = {10.1103/PhysRevE.66.046129},
  url = {https://link.aps.org/doi/10.1103/PhysRevE.66.046129}
}

\bib{MR3796372}{article}{
   author={Cusimano, Nicole},
   author={del Teso, F\'{e}lix},
   author={Gerardo-Giorda, Luca},
   author={Pagnini, Gianni},
   title={Discretizations of the spectral fractional Laplacian on general
   domains with Dirichlet, Neumann, and Robin boundary conditions},
   journal={SIAM J. Numer. Anal.},
   volume={56},
   date={2018},
   number={3},
   pages={1243--1272},
   issn={0036-1429},
   review={\MR{3796372}},
   doi={10.1137/17M1128010},
}

\bib{VERO}{article}{
   author={Dipierro, Serena},
   author={Proietti Lippi, Edoardo},
   author={Valdinoci, Enrico},
title={(Non)local logistic equations with Neumann conditions},
journal = {arXiv e-prints},
     date = {2021},
          eid = {arXiv:2101.02315},
        pages = {arXiv:2101.02315},
archivePrefix = {arXiv},
      adsnote = {Provided by the SAO/NASA Astrophysics Data System}  
      }

\bib{dirova}{article}{
   author={Dipierro, Serena},
   author={Ros-Oton, Xavier},
   author={Valdinoci, Enrico},
   title={Nonlocal problems with Neumann boundary conditions},
   journal={Rev. Mat. Iberoam.},
   volume={33},
   date={2017},
   number={2},
   pages={377--416},
   issn={0213-2230},
   doi={10.4171/RMI/942},
}

\bib{ELEME}{article}{
author = {Dipierro, Serena},
author = {Valdinoci, Enrico},
        title = {Elliptic partial differential equations from an elementary viewpoint},
      journal = {arXiv e-prints},
date = {2020},
          eid = {arXiv:2101.07941},
        pages = {arXiv:2101.07941},
archivePrefix = {arXiv},
      adsnote = {Provided by the SAO/NASA Astrophysics Data System}
}

\bib{NATUBIS}{article}{
   author={Edwards, A. M.},
   author={Phillips, R. A.},
      author={Watkins, N. W.},
         author={Freeman, M. P.},
            author={Murphy, E. J.},
               author={Afanasyev, V.},
                  author={Buldyrev, S. V.},
                     author={da Luz, M. G. E.},
                        author={Raposo, E. P.},
                           author={Stanley, H. E.},
                              author={Viswanathan, G. M.},
   title={Revisiting L\'evy flight search patterns of wandering albatrosses, bumblebees and deer},
   journal={Nature},
   volume={449},
   date={2007},
   pages={1--34},
   issn={1476-4687},
   doi={10.1038/nature06199},
   }

\bib{MR2743162}{book}{
   author={Fischer, Hans},
   title={A history of the central limit theorem},
   series={Sources and Studies in the History of Mathematics and Physical
   Sciences},
   note={From classical to modern probability theory},
   publisher={Springer, New York},
   date={2011},
   pages={xvi+402},
   isbn={978-0-387-87856-0},
   doi={10.1007/978-0-387-87857-7},
}

\bib{MR0233400}{book}{
   author={Gnedenko, B. V.},
   author={Kolmogorov, A. N.},
   title={Limit distributions for sums of independent random variables},
   series={Translated from the Russian, annotated, and revised by K. L.
   Chung. With appendices by J. L. Doob and P. L. Hsu. Revised edition},
   publisher={Addison-Wesley Publishing Co., Reading, Mass.-London-Don
   Mills., Ont. },
   date={1968},
   pages={ix+293},
}

\bib{SHARK}{article}{
author={Hale, Tom},
title={Great white shark dies after just three days in captivity},
year={2016},
journal={IFL Science,
https://www.iflscience.com/plants-and-animals/great-white-shark-dies-after-just-three-days-captivity/},
}

\bib{FALCO}{article}{
author={Hern\'andez-Pliego, Jes\'us},
author={Rodr\'{\i}guez, Carlos},
author={Bustamante, Javier},
date={2017},
title={A few long versus many short foraging trips: different foraging strategies of lesser kestrel sexes during breeding}
journal={Movement Ecology},
pages={1--11},
number={7},
issue={41},
ISSN={2051-3933},
url={https://doi.org/10.1186/s40462-017-0100-6},
doi={10.1186/s40462-017-0100-6},
}

\bib{MR837810}{article}{
   author={Hsu, Pei},
   title={On excursions of reflecting Brownian motion},
   journal={Trans. Amer. Math. Soc.},
   volume={296},
   date={1986},
   number={1},
   pages={239--264},
   issn={0002-9947},
   doi={10.2307/2000572},
}

\bib{NATU2}{article}{
   author={Humphries, Nicolas E.},
   author={Queiroz, Nuno},
   author={Dyer, Jennifer R. M.},
   author={Pade, Nicolas G.},
   author={Musyl, Michael K.},
   author={Schaefer, Kurt M.},
   author={Fuller, Daniel W.},
   author={Brunnschweiler, Juerg M.},
   author={Doyle, Thomas K.},
   author={Houghton, Jonathan D. R.},
   author={Hays, Graeme C.},
   author={Jones, Catherine S.},
   author={Noble, Leslie R.},
   author={Wearmouth, Victoria J.},
   author={Southall, Emily J.},
   author={Sims, David W.},
   title={Environmental context explains L\'evy and Brownian movement patterns of marine predators},
   journal={Nature},
   volume={465},
   date={2010},
   pages={1066--1069},
   issn={1476-4687},
   doi={10.1038/nature09116},
}

\bib{PNASS2}{article}{
	author = {Humphries, Nicolas E.},
		author = {Weimerskirch, Henri},
			author = {Queiroz, Nuno},
				author = {Southall, Emily J.},
					author = {Sims, David W.},
	title = {Foraging success of biological L{\'e}vy flights recorded in situ},
	volume = {109},
	number = {19},
	pages = {7169--7174},
	year = {2012},
	doi = {10.1073/pnas.1121201109},
	publisher = {National Academy of Sciences},
	issn = {0027-8424},
	URL = {https://www.pnas.org/content/109/19/7169},
	eprint = {https://www.pnas.org/content/109/19/7169.full.pdf},
	journal = {Proceedings of the National Academy of Sciences}
}

\bib{MR2609393}{article}{
   author={James, Alex},
   author={Pitchford, Jonathan W.},
   author={Plank, Michael J.},
   title={Efficient or inaccurate? Analytical and numerical modelling of
   random search strategies},
   journal={Bull. Math. Biol.},
   volume={72},
   date={2010},
   number={4},
   pages={896--913},
   issn={0092-8240},
   doi={10.1007/s11538-009-9473-z},
}

\bib{PhysRevE78051128}{article}{
  title = {Optimizing the encounter rate in biological interactions: Ballistic versus L\'evy versus Brownian strategies},
  author = {James, A.},
    author = {Plank, M. J.},
      author = {Brown, R.},
  journal = {Phys. Rev. E},
  volume = {78},
  issue = {5},
  pages = {051128},
  numpages = {5},
  year = {2008},
  doi = {10.1103/PhysRevE.78.051128},
  url = {https://link.aps.org/doi/10.1103/PhysRevE.78.051128}
}

\bib{Klages R}{article}{ author={Klages, Rainer}, title={Search for Food of Birds, Fish and Insects}, conference={ title={Diffusive Spreading in Nature, Technology and Society}, }, book={ publisher={Springer, Cham}, }, date={2018}, pages={1--21}, doi={10.1007/978-3-319-67798-9\_4},}

\bib{MR3613319}{article}{
   author={Kwa\'{s}nicki, Mateusz},
   title={Ten equivalent definitions of the fractional Laplace operator},
   journal={Fract. Calc. Appl. Anal.},
   volume={20},
   date={2017},
   number={1},
   pages={7--51},
   issn={1311-0454},
   doi={10.1515/fca-2017-0002},
}

\bib{FILO}{article}{
    author = {Landis, Michael J.},
    author = {Schraiber, Joshua G.},
    author = {Liang, Mason},
    title = {Phylogenetic analysis using L\'evy Processes: finding jumps in the evolution of continuous traits},
    journal = {Systematic Biology},
    volume = {62},
    number = {2},
    pages = {193--204},
    year = {2012},
    issn = {1063-5157},
    doi = {10.1093/sysbio/sys086},
    url = {https://doi.org/10.1093/sysbio/sys086},
}

\bib{PhysRevLett108098103}{article}{
author = {Lenz, Friedrich},
author = {Ings, Thomas C.},
author = {Chittka, Lars},
author = {Chechkin, Aleksei V.},
author = {Klages, Rainer},
  title = {Spatiotemporal Dynamics of Bumblebees Foraging under Predation Risk},
  journal = {Phys. Rev. Lett.},
  volume = {108},
  issue = {9},
  pages = {098103},
  numpages = {5},
  year = {2012},
  doi = {10.1103/PhysRevLett.108.098103},
  url = {https://link.aps.org/doi/10.1103/PhysRevLett.108.098103}
}

\bib{LIT1}{article}{
  title = {Inverse Square L\'evy Walks are not Optimal Search Strategies for $d{\ge}2$},
  author = {Levernier, Nicolas},
    author = {Textor, Johannes},
      author = {B\'enichou, Olivier},
        author = {Voituriez, Rapha\"el},
  journal = {Phys. Rev. Lett.},
  volume = {124},
  issue = {8},
  pages = {080601},
  numpages = {5},
  year = {2020},
  doi = {10.1103/PhysRevLett.124.080601},
  url = {https://link.aps.org/doi/10.1103/PhysRevLett.124.080601}
}

\bib{LIT3}{article}{
  title = {Reply to ``Comment on `Inverse Square L\'evy Walks are not Optimal Search Strategies for $d{\ge}2$'''},
  author = {Levernier, Nicolas},
    author = {Textor, Johannes},
      author = {B\'enichou, Olivier},
        author = {Voituriez, Rapha\"el},
  journal = {Phys. Rev. Lett.},
  volume = {126},
  issue = {4},
  pages = {048902},
  numpages = {1},
  year = {2021},
  doi = {10.1103/PhysRevLett.126.048902},
  url = {https://link.aps.org/doi/10.1103/PhysRevLett.126.048902}
}

\bib{MR2459736}{article}{
   author={Mainardi, Francesco},
   author={Pagnini, Gianni},
   title={Mellin-Barnes integrals for stable distributions and their
   convolutions},
   journal={Fract. Calc. Appl. Anal.},
   volume={11},
   date={2008},
   number={4},
   pages={443--456},
   issn={1311-0454},
   review={\MR{2459736}},
}

\bib{MR3246044}{article}{
   author={Musina, Roberta},
   author={Nazarov, Alexander I.},
   title={On fractional Laplacians},
   journal={Comm. Partial Differential Equations},
   volume={39},
   date={2014},
   number={9},
   pages={1780--1790},
   issn={0360-5302},
   doi={10.1080/03605302.2013.864304},
}

\bib{MR3886705}{article}{
   author={Nadarajah, Saralees},
   author={Chan, Stephen},
   title={The exact distribution of the sum of stable random variables},
   journal={J. Comput. Appl. Math.},
   volume={349},
   date={2019},
   pages={187--196},
   issn={0377-0427},
   review={\MR{3886705}},
   doi={10.1016/j.cam.2018.09.044},
}

\bib{MR2559350}{article}{
   author={Pagnini, Gianni},
   author={Mainardi, Francesco},
   title={Evolution equations for the probabilistic generalization of the
   Voigt profile function},
   journal={J. Comput. Appl. Math.},
   volume={233},
   date={2010},
   number={6},
   pages={1590--1595},
   issn={0377-0427},
   review={\MR{2559350}},
   doi={10.1016/j.cam.2008.04.040},
}

\bib{Palyulin2931}{article}{
	author = {Palyulin, Vladimir V.},
		author = {Chechkin, Aleksei V.},
			author = {Metzler, Ralf},
	title = {L{\'e}vy flights do not always optimize random blind search for sparse targets},
	volume = {111},
	number = {8},
	pages = {2931--2936},
	year = {2014},
	doi = {10.1073/pnas.1320424111},
	issn = {0027-8424},
	URL = {https://www.pnas.org/content/111/8/2931},
	eprint = {https://www.pnas.org/content/111/8/2931.full.pdf},
	journal = {Proceedings of the National Academy of Sciences}
}

\bib{POL}{article}{
author = {Polverino, Giovanni},
author = {Manciocco, Arianna},
author={Alleva, Enrico},
title = {Effects of spatial and social restrictions on the presence of stereotypies in the budgerigar (Melopsittacus undulatus): a pilot study},
journal = {Ethology Ecology \& Evolution},
volume = {24},
number = {1},
pages = {39--53},
year  = {2012},
doi = {10.1080/03949370.2011.582045},
       adsurl ={https://doi.org/10.1080/03949370.2011.582045}    
}

\bib{POL2}{article}{
author = {Polverino, Giovanni},
author = {Manciocco, Arianna},
author={Vitale, Augusto},
author={Alleva, Enrico},
title = {Stereotypic behaviours in Melopsittacus undulatus: Behavioural consequences of social and spatial limitations},
journal = {Applied Animal Behaviour Science},
volume = {165},
pages = {143--155},
year = {2015},
issn = {0168-1591},
doi = {https://doi.org/10.1016/j.applanim.2015.02.009},
adsurl ={http://www.sciencedirect.com/science/article/pii/S016815911500057X},
}

\bib{MR2639124}{article}{
   author={Revelli, J. A.},
   author={Rojo, F.},
   author={Budde, C. E.},
   author={Wio, H. S.},
   title={Optimal intermittent search strategies: {\it smelling} the prey},
   journal={J. Phys. A},
   volume={43},
   date={2010},
   number={19},
   pages={195001, 11},
   issn={1751-8113},
   doi={10.1088/1751-8113/43/19/195001},
}

\bib{MR2670512}{article}{
   author={Rojo, F.},
   author={Revelli, J.},
   author={Budde, C. E.},
   author={Wio, H. S.},
   author={Oshanin, G.},
   author={Lindenberg, Katja},
   title={Intermittent search strategies revisited: effect of the jump
   length and biased motion},
   journal={J. Phys. A},
   volume={43},
   date={2010},
   number={34},
   pages={345001, 10},
   issn={1751-8113},
   doi={10.1088/1751-8113/43/34/345001},
}

\bib{ANT}{article}{
author={Saar, Maya}
author={Subach, Aziz},
author={Reato, Illan},
author={Liber, Tal}
author={Pruitt, Jonathan N.},
author={Scharf, Inon},
title={Consistent differences in foraging behavior in 2 sympatric harvester ant species may facilitate coexistence},
issn={1674-5507},
date={2018}
pages={653--661},
volume={64},
issue={5},
URL={https://pubmed.ncbi.nlm.nih.gov/30323844},
DOI={10.1093/cz/zox054},
journal={Current Zoology},
}

\bib{MR2399851}{book}{
   author={Salsa, Sandro},
   title={Partial differential equations in action},
   series={Universitext},
   note={From modelling to theory},
   publisher={Springer-Verlag Italia, Milan},
   date={2008},
   pages={xvi+556},
   isbn={978-88-470-0751-2},
}

\bib{MR3380662}{book}{
   author={Salsa, Sandro},
   author={Verzini, Gianmaria},
   title={Partial differential equations in action},
   series={Unitext},
   volume={87},
   note={Complements and exercises;
   Translated from the 2005 Italian edition by Simon G. Chiossi;
   La Matematica per il 3+2},
   publisher={Springer, Cham},
   date={2015},
   pages={viii+422},
   isbn={978-3-319-15415-2},
   isbn={978-3-319-15416-9},
   doi={10.1007/978-3-319-15416-9},
}

\bib{MR3233760}{article}{
   author={Servadei, Raffaella},
   author={Valdinoci, Enrico},
   title={On the spectrum of two different fractional operators},
   journal={Proc. Roy. Soc. Edinburgh Sect. A},
   volume={144},
   date={2014},
   number={4},
   pages={831--855},
   issn={0308-2105},
   doi={10.1017/S0308210512001783},
}

\bib{SIMS}{article}{
author = {Sims, David W.},
author = {Humphries, Nicolas E.},
author = {Bradford, Russell W.},
author = {Bruce, Barry D.},
title = {L\'evy flight and Brownian search patterns of a free-ranging predator reflect different prey field characteristics},
journal = {Journal of Animal Ecology},
volume = {81},
number = {2},
pages = {432-442},
doi = {https://doi.org/10.1111/j.1365-2656.2011.01914.x},
url = {https://besjournals.onlinelibrary.wiley.com/doi/abs/10.1111/j.1365-2656.2011.01914.x},
eprint = {https://besjournals.onlinelibrary.wiley.com/doi/pdf/10.1111/j.1365-2656.2011.01914.x},
year = {2012}
}

\bib{SKORO1}{article}{
author={Skorokhod, Anatoliy V.},
date={1961},
title={Stochastic equations for diffusion processes in a bounded region 1},
journal={Theor. Veroyatnost. i Primenen}
volume={6},
pages={264--274},
}

\bib{SKORO2}{article}{
author={Skorokhod, Anatoliy V.},
date={1962},
title={Stochastic equations for diffusion processes in a bounded region 2},
journal={Theor. Veroyatnost. i Primenen}
volume={7},
pages={3--23.},
}

\bib{MR3590646}{article}{
   author={Sprekels, J\"{u}rgen},
   author={Valdinoci, Enrico},
   title={A new type of identification problems: optimizing the fractional
   order in a nonlocal evolution equation},
   journal={SIAM J. Control Optim.},
   volume={55},
   date={2017},
   number={1},
   pages={70--93},
   issn={0363-0129},
   review={\MR{3590646}},
   doi={10.1137/16M105575X},
}

\bib{MR2584076}{article}{
   author={Valdinoci, Enrico},
   title={From the long jump random walk to the fractional Laplacian},
   journal={Bol. Soc. Esp. Mat. Apl. SeMA},
   number={49},
   date={2009},
   pages={33--44},
   issn={1575-9822},
}

\bib{VISWANATHAN2001}{article}{
   title = {Statistical physics of random searches},
   journal = {Brazilian Journal of Physics},
   author={Viswanathan, G. M.},
      author={Afanasyev, V.},
         author={Buldyrev, Sergey V.},
            author={Havlin, Shlomo},
               author={da Luz, M. G. E.},
                  author={Raposo, E. P.},
                     author={Stanley, H. Eugene},
   ISSN = {0103-9733},
   URL = {http://www.scielo.br/scielo.php?script=sci_arttext&pid=S0103-97332001000100018&nrm=iso},
   volume = {31},
   year = {2001},
   pages = {102 - 108},
   doi = {10.1590/S0103-97332001000100018},
   }

\bib{NATU}{article}{
   author={Viswanathan, G. M.},
   author={Afanasyev, V.},
   author={Buldyrev, S. V.},
   author={Murphy, E. J.},
   author={Prince, P. A.},
   author={Stanley, H. E.},
   title={L\'evy flight search patterns of wandering albatrosses},
   journal={Nature},
   volume={381},
   date={1996},
   pages={413--415},
   issn={1476-4687},
   doi={10.1038/381413a0},
}

\bib{89897765NATU}{article}{
   author={Viswanathan, G. M.},
   author={Buldyrev, S. V.},
   author={Havlin, S.},
   author={da Luz, M. G. E.},
   author={Raposo, E. P.},
   author={Stanley, H. E.},
   title={Optimizing the success of random searches},
   journal={Nature},
   volume={401},
   date={1999},
   pages={911--914},
   issn={1476-4687},
   doi={10.1038/44831},
}

\bib{FORA}{book}{
publisher={Cambridge University Press, Cambridge},
title={The Physics of Foraging: An Introduction to Random Searches and Biological Encounters}, 
DOI={10.1017/CBO9780511902680}, 
author={Viswanathan, Gandhimohan. M.},
author={da Luz, Marcos G. E.},
author={Raposo, Ernesto P.},
author={Stanley, H. Eugene},
isbn={9780511902680}
date={2011},
   pages={xiv+164},
}

\bib{VOND}{article}{
       author = {Vondra{\v{c}}ek, Zoran},
        title = {A probabilistic approach to a non-local quadratic
        form and its connection to the Neumann boundary condition problem},
      journal = {arXiv e-prints},
     date = {2019},
          eid = {arXiv:1909.10687},
        pages = {arXiv:1909.10687},
archivePrefix = {arXiv},
      adsnote = {Provided by the SAO/NASA Astrophysics Data System}
}

\bib{MR3403266}{article}{
   author={Zaburdaev, V.},
   author={Denisov, S.},
   author={Klafter, J.},
   title={L\'{e}vy walks},
   journal={Rev. Modern Phys.},
   volume={87},
   date={2015},
   number={2},
   pages={483--530},
   issn={0034-6861},
   doi={10.1103/RevModPhys.87.483},
}

\end{biblist}
\end{bibdiv}
\end{document}